\documentclass[10pt]{article}

\newcommand{\mylbl}[1]{{\label{#1}}}

\newcommand{\mbold}[1]{\mbox{\bf #1}}
\newcommand{\C}{{\mbold C}}

\newcommand{\N}{{\mbold N}}
\newcommand{\Z}{{\mbold Z}}
\newcommand{\hf}[1]{{\frac{#1}{2}}}
\newcommand{\bdot}{{\,\stackrel{\scriptstyle\bullet}{{}}\,}}
\newcommand{\mdot}{{\,\stackrel{\scriptstyle\circ}{{}}\,}}
\newcommand{\bysame}{\leavevmode\hbox to 3em{\hrulefill}}

\newcommand{\CP}{{\cal P}}

\newcommand{\Cg}{{\cal G}}
\newcommand{\CF}{{\cal F}}
\newcommand{\CW}{{\cal W}}

\newcommand{\RC}{{\check R}}
\newcommand{\IC}{{\check I}{}}
\newcommand{\JC}{{\check J}{}}
\newcommand{\VC}{P}
\newcommand{\RH}{{\hat R}}

\newcommand{\Span}{{\mbox{\rm Span}}}
\newcommand{\Res}{{\mbox{\rm Res}}}
\newcommand{\Cl}{{\mbox{{\rm Cl}}}}
\newcommand{\wt}{{{\Delta}}}

\newcommand{\End}{{\mbox{\rm End}}}
\newcommand{\Hom}{{\mbox{\rm Hom}}}

\newcommand{\Img}{{\mbox{\rm Im}}}
\newcommand{\Vir}{{\mbox{\it Vir}}}
\newcommand{\CKS}{{{C\!K_6}}}
\newcommand{\NFA}{{{N_4^\alpha}}}
\newcommand{\NFB}{{{N_4^\beta}}}
\newcommand{\K}{\C}
\newcommand{\ba}{{{\alpha}}}

\newcommand{\Eb}[2]{{\bar{D}_{#1}\mdot \bar{D}_{#2}}}
\newcommand{\F}[2]{{\bar{D}_{#1}\mdot D_{#2}}}

\newcommand{\ZZ}{{\Z/2\Z}}
\newcommand{\myqed}{\hfill \hbox{\rule[-2pt]{3pt}{6pt}}\newline}

\newcommand{\pprod}[1]{{\langle{#1}\rangle}}
\newcommand{\pp}[3]{{{#1}{}_{\pprod{#2}}{#3}}}
\newcommand{\ppp}[3]{{{#1}{}_{\pprod{#2}}{#3}}}
\newcommand{\PP}[1]{{\pprod{#1}}}
\newcommand{\PL}[1]{{#1}/{\partial{#1}}}

\newtheorem{definition}{Definition}
\newtheorem{proposition}{Proposition}
\newtheorem{theorem}{Theorem}
\newtheorem{example}{Example}
\newtheorem{remark}{Remark}
\newtheorem{lemma}{Lemma}
\newtheorem{note}{Note}
\newtheorem{corollary}{Corollary}

\newcommand{\pr}{{\sc Proof} \quad}

\newenvironment{namelist}[1]{\begin{list}{}{\settowidth{\labelwidth}{#1}\setlength{\leftmargin}{1.1\labelwidth}}}{\end{list}}

\makeatletter
\newcounter{mypcount}
 \@addtoreset{mypcount}{section}
\makeatother
\makeatletter
  
  \@addtoreset{equation}{section}
\makeatother
\makeatletter
  
  \@addtoreset{theorem}{section}
\makeatother
\makeatletter
  
  \@addtoreset{proposition}{section}
\makeatother
\makeatletter
  
  \@addtoreset{lemma}{section}
\makeatother
\makeatletter
  
  \@addtoreset{corollary}{section}
\makeatother
\makeatletter
  
  \@addtoreset{definition}{section}
\makeatother

\makeatletter
  
  \@addtoreset{remark}{section}
\makeatother
\makeatletter
  
  \@addtoreset{note}{section}
\makeatother
\makeatletter
  
  \@addtoreset{example}{section}
\makeatother

\begin{document}
\title{\large\bf Algebraic structures on quasi-primary states in superconformal algebras}
\author{Go Yamamoto\footnote{Email: yamamo@ms.u-tokyo.ac.jp}\\ \\{\it
Graduate School of Mathematical Sciences, The University of Tokyo}\\{\it 3-8-1 Komaba, Meguro-ku, Tokyo 153-8914, Japan.}}
\date{}
\maketitle
\begin{abstract}
The algebraic structure on the subspace of the quasi-primary vectors given by the projection of the $(n)$ products of a conformal superalgebra is formulated.
 As an application the complete list of simple physical conformal superalgebras is given.  The list contains a one-parameter family of superconformal algebras with $4$ supercharges that is simple for general values.
\end{abstract}

\section{Introduction}

  For an infinite-dimensional Lie superalgebra $\Cg$, one often assumes that there exists a finite set $\CF$ of generating functions of elements of $\Cg$ and that the Lie bracket is written in terms of the OPE (Operator Product Expansion), i.e.,
\begin{equation}
a(z)b(w) \sim \sum_j\frac{c^j(z)}{(z-w)^j},
\end{equation}
where the $\sum$ is finite.  It means 
\begin{eqnarray}
 [a(z),b(w)] = \sum_j\frac{c^j(w)}{j!}\partial^j\delta(z-w),\\
 c^j(w)= \Res_z[a(z),b(w)](z-w)^j,
\end{eqnarray}
for $a,b \in \C[\partial]\CF$, where the $\sum$ is always finite.  The finiteness is called locality.  Many significant infinite-dimensional Lie superalgebras, e.g., affine Lie algebras, the Virasoro algebra, the Neveu-Schwarz algebra, have locality.

  The notion of conformal superalgebra (vertex Lie superalgebra) is formulated in \cite{kac1} and \cite{Primc} independently, which is an axiomatic description of Lie superalgebras with OPE with respect to the infinitely many operations $a_{(j)}b= c^j$ as above.  Once a conformal superalgebra is given, one can reconstruct the Lie superalgebra $\Cg$.  We shall require existence of conformal vector in addition, which corresponds to a Virasoro subalgebra in the associated Lie superalgebra.

 For a conformal superalgebra $R$ the subspace of the quasi-primary vectors (see section 2 for precise formulation) are identified with $R/(\partial R)$.  
  For some kind of conformal superalgebras the space of primary states generates the conformal superalgebra and the associated Lie superalgebra (\cite{Primc}).
  On the other hand, it is well-known that  a conformal superalgebra $R$ yields a Lie superalgebra structure on $R/(\partial R)$.
  We will study more detailed structures on $R/(\partial R)$ (section 3).

The algebraic structures on the space of the primary vectors are described in \cite{BPZ}.   We will study the algebraic structures on the space of the quasi-primary vectors, defining the $\PP{n}$ products on it by the projection of the $(n)$ products.   We will show that one can reconstruct the entire conformal superalgebra from the $\PP{n}$ products on the space of the quasi-primary vectors (section 4).  For the normal product the structure given by the projection are discussed in \cite{proj}.

  The $\PP{n}$ products are (anti-)commutative, but may not be associative.  We have some simple cases of the algebraic structures on the subspace of the quasi-primary vectors.  The most simple one is the case of affine Lie superalgebras, for which all products but the $\PP{0}$ product vanish and the $\PP{0}$ product yields a Lie superalgebra structure on the finite-dimensional vector space of the quasi-primary vectors.
  The second case is physical conformal superalgebras, which corresponds to the superconformal algebras, for example, the Virasoro algebra, the Neveu-Schwarz algebra and the $N=4$ superconformal algebra.  In this case all products but the $\PP{0}$ product and the $\PP{1}$ product vanish and they
 yield a left Clifford module structure on the finite-dimensional vector space of the quasi-primary vectors.  
The action of Clifford algebra is described in \cite{HTT}.  It restricts the dimension of the space of the quasi-primary vectors of physical conformal superalgebras.

Examples of simple physical conformal superalgebras are given in \cite{ck}, \cite{KL}, \cite{kac1}, and \cite{kac2}.  The list of known simple physical conformal superalgebras are $\Vir$, $K_1$, $K_2$, $K_3$, $S_2$, $W_2$, $\CKS$, where we have followed the notations of \cite{ck} and \cite{kac0}.  $\Vir$ is the Virasoro algebra.  $K_j$ is known as the $N=j$ superconformal algebra.  $S_2$ and $W_2$ are superconformal algebras with $4$ supercharges.  $S_2$ is known as the $N=4$ superconformal algebra.  $\CKS$ is discovered in \cite{ck} and is the only known superconformal algebra with more than $4$ supercharges.

In \cite{kac2} a list for the simple physical conformal superalgebras is given, however, we are making another approach.
As an application of the reconstruction theorem we will classify simple physical conformal superalgebras by working on the space of the quasi-primary vectors and the $\PP{n}$ products on it.  We have found a simple physical conformal superalgebra $N_4$ and a one-parameter family of physical conformal superalgebras $\NFA$ that is simple for all $\alpha\in(\C/\{\pm1\})\setminus\{[1]\}$, which imply a class of simple physical conformal algebras that is not in the list of \cite{kac2} exists; $N_4$ and $\NFA$s are counter examples to Lemma 4.1(b) in \cite{kac2}.  The simple physical conformal superalgebras $N_4$ and $\NFA$ coincide with the centerless conformal superalgebras of the large $N=4$ superconformal algebras written down in \cite{STP}.
The complete list of the simple physical conformal superalgebras is $\Vir$, $K_1$, $K_2$, $K_3$, $S_2$, $W_2$, $N_4$, $\NFA$ and $\CKS$ (section 8).

 {\bf Acknowledgements}
I deeply wish to thank Atsushi Matsuo, my research supervisor, for suggesting many of the problems and for encouraging me to keep the directions pursued in this paper.
I am very grateful to Shun-Jen Cheng for valuable comments on changing the conformal vectors of $\NFA$.
I also wish to thank Kenji Iohara, Masao Jinzenji, Akishi Kato and Minoru Wakimoto for valuable discussions and comments.

\section{ Preliminaries}

 Let $K$ be a subfield of $\C$. 
A $K$-vector space $V$ with a direct sum decomposition $V=V_0 \oplus V_1$ is called a $\ZZ$-graded $K$-vector space. 
The homomorphisms of $\ZZ$-graded vector spaces are supposed to be compatible with the gradation.
 The $\ZZ$-gradation is called {\it parity}. 
$V_0$ is called the subspace of {\it even} parity, and $V_1$ is of {\it odd} parity. 

The $\ZZ$-graded objects are called {\it super-} objects.
  Commutativity for the product $\cdot$ of a superalgebra is defined to be $ a \cdot b = (-1)^{p(a)p(b)} b \cdot a ,$ where $a,b$ are supposed to be homogeneous with respect to the parity $p$.

Now let us state the axioms for conformal superalgebras, based on the descriptions in \cite{kac1} and \cite{kac2}.
We denote $A^{(j)}={A^j}/{j!}$, where $A$ is an operator.

\addtocounter{mypcount}{1}\begin{definition}{\rm \mylbl{basics-df-1}
 Let $R$ be a $\ZZ$-graded $K$-vector space equipped with countably many products
$$(n) : R \otimes R \to R, \quad (n\in\N),  $$
and a linear map $\partial : R \to R$.
 The triple $(R,\{(n)\}_{n\in\N},L)$ satisfying the following conditions for an even vector $L\in R$ are called a {\it conformal superalgebra}:
\begin{namelist}{(Cx)}
\item[(C)] For all $a,b,c \in R$,
\item[(C0)] there exists some $N \in \N$ such that for all $ n \in \N$ satisfying $n\ge N$
$$a_{(n)}b=0  ,$$
\item[(C1)] for all $n \in \N$,
 $$ (\partial a)_{(n)}b=-na_{(n-1)}b  ,$$ 
\item[(C2)] for all $n \in \N$,
 $$ a_{(n)}b=(-1)^{p(a)p(b)}\sum\limits_{j=0}^{\infty} {(-1)^{j+n+1}\partial ^{(j)}b_{(n+j)}a}  ,$$
\item[(C3)] for all $m,n \in \N$,
 $$  a_{(m)}(b_{(n)}c)=\sum\limits_{j=0}^{\infty} {\left( {\matrix{m\cr
j\cr
}} \right)\left( {a_{(j)}b} \right)_{(n+m-j)}c}+(-1)^{p(a)p(b)}b_{(n)}(a_{(m)}c)  .$$
\item[(V)] $L \in R$ satisfies
$ L_{(0)}L = \partial L$,
$ L_{(1)}L = 2L$,
$ L_{(2)}L = 0$,
$L_{(0)} = \partial$ as operators on $R$,
and $L_{(1)}$ is diagonalizable.
\end{namelist}
}\end{definition}

\addtocounter{mypcount}{1}\begin{remark}{\rm \mylbl{basics-rem-2}
The $\sum$ in (C3) is a finite sum because of (C0).
}\end{remark}

$L$ is called the {\it conformal vector} of $R$.
A homomorphism of conformal superalgebras from  $R$ to $R'$ is a $K[\partial]$-module homomorphism  $f: R \to R'$  that is compatible with the $(n)$ products for all $n \in \N$ and maps $L$ to the conformal vector of $R'$.
 An ideal of a conformal superalgebra is a $K[\partial]$-submodule that is closed under the left multiplication of the $(n)$ products for all $n \in \N$.
 A conformal superalgebra $R$ with no ideals other than $\{0\}$ and $R$ itself is called a {\it simple} conformal superalgebra.
 The ideal  $ \{ c \in R | x_{(n)}c=0, x \in R, n \in \N\} $ is called the {\it center} of $ R $.
 If the center is $ \{0\}$  then the conformal superalgebra is said to be {\it centerless}.

\addtocounter{mypcount}{1}\begin{remark}{\rm \mylbl{basics-rem-5}
Right ideals are defined similarly, but they coincides with left ideals.
}\end{remark}

\addtocounter{mypcount}{1}\begin{note}{\rm
The axiom (V) is not included in the definition of conformal superalgebras in \cite{kac1} and \cite{kac2} nor of vertex Lie algebras in \cite{Primc}, while existence of the conformal vector is assumed for superconformal algebras.  We set it into the axioms for conformal superalgebras.
}\end{note}

Suppose given an isomorphism of $K[\partial]$-modules $f: R \to R'$ that is compatible with the $(n)$ products where $(R,\{(n)\}_{n\in\N},L)$ and $(R',\{(n)\}_{n\in\N},L')$ are conformal superalgebras.  We say that $(R,\{(n)\}_{n\in\N},L)$ and $(R',\{(n)\}_{n\in\N},L')$ are {\it equivalent} to each other.

\addtocounter{mypcount}{1}\begin{note}{\rm \mylbl{basics-rem-equiv}
In \cite{kac1} and \cite{kac2} the class with respect to the equivalence above is considered.  We will consider the isomorphism classes, which is stronger than to the equivalence classes.
}\end{note}

The eigenvalue of $L_{(1)}$ is denoted by $\wt(x)$ for an eigenvector $x$ and is called the {\it conformal weight of $x$}. 
Define
$ R^k = \{ x \in R | L_{(1)}x = k x\}$, $\wt_R = \{ k \in K | R^k \ne \{0\}\}$ and $\wt'_R=\wt_R\setminus\{0\}$.

\addtocounter{mypcount}{1}\begin{remark}{\rm \mylbl{basics-rem-16}
We have
$\wt(\partial x) = \wt(x) +1$ and
$\wt(x_{(n)}y) = \wt(x) + \wt(y) -n -1$.
That is,\\ $(n): R^p \otimes R^q \to R^{p+q-n-1}$.
}\end{remark}

 A conformal superalgebra $R$ over $\C$ is called a {\it superconformal algebra} if there exists a finite-dimensional subspace $\CF$ such that $R =  \C[\partial] \CF$, all conformal weights are non-negative half-integers, the even subspace $R_{\mbox{even}} = \bigoplus_{n\in\N}R^n$ and the odd subspace $R_{\mbox{odd}} = \bigoplus_{n\in\N+\hf{1}}R^n$.
 We call a superconformal algebra $R$ a {\it physical conformal superalgebra} if $\CF \subset R^2 \oplus R^\hf{3} \oplus R^1 \oplus R^\hf{1}$ and $\CF \cap R^2 = \C L$, following the terminology in \cite{kac2}.

Set $ \RH = \bigoplus_{n\in\Z} R_{(n)} $ and $J = \Span\{ (\partial a)_{(n)} + n a_{(n-1)}|\ a \in R,\ n \in \Z\}$ where $ R_{(n)} $ is a copy of $ R $ for each $n\in\Z$.
The Lie superalgebra $ \RH /J $ defined by
 $ [a_{(m)},b_{(n)}]=\sum_{j=0}^{\infty} {\left( {\matrix{m\cr
j\cr
}} \right)\left( {a_{(j)}b} \right)_{(n+m-j)}} $
 is called {\it the Lie superalgebra associated to} $R$.  If the conformal superalgebra $R$ is not simple then
the Lie superalgebra associated to $R$ is not simple.

\section{ $\partial$-decomposition }

For a conformal superalgebra $(R,\{(n)\}_{n\in\N},L)$, we shall call the subspace 
$ \{ {x \in R} | {L_{(2)}x} \in {R^0} \} $
the {\it reduced subspace} of $R$ and denote it by $\RC$.  We call the elements of the reduced subspace {\it reduced vectors}.  Denote $\RC^k = \RC \cap R^k$, $\wt_\RC=\{k\in K| \RC^k \ne \{0\}\}$ and $\wt'_\RC=\wt_\RC\setminus\{0\}$.  Obviously we have $\RC^0 = R^0$.

\addtocounter{mypcount}{1}\begin{remark}{\rm \mylbl{decomp-rem-15}
If the Lie algebra $(R^1,(0))$ is perfect then we have $L_{(2)}a = 0$ for all $a\in \RC$ because $L_{(2)} a_{(0)}b = 0$ for all $a,b \in \RC^1$.
}\end{remark}

We introduce the notion of {\it regular} conformal superalgebras.  The superconformal algebras are regular.

\addtocounter{mypcount}{1}\begin{definition}{\rm \mylbl{basics-df-11}
A conformal superalgebra $R$ is regular if $R^0$ is the center and if $\wt_R\cap(-\hf{1}\N) \subset \{0\}$ and for each $k \in \wt_R$ there exists some $M \in\N$ such that $k-m \notin \wt_R$ for all $m\in\N$ satisfying $m\ge M$.
}\end{definition}

\addtocounter{mypcount}{1}\begin{proposition}\mylbl{decomp-prop-1}
Let $(R,\{(n)\}_{n\in\N},L)$ be a regular conformal superalgebra and $\RC$ the reduced subspace of $(R,\{(n)\}_{n\in\N},L)$.
Then
there exists a unique decomposition
\begin{equation}
 x = \sum_{j=0}^{m} \partial^{(j)} x^j 
\end{equation}
for any $x \in R$ for some $m\in\N$ where  $x^0 \in \RC$ and $x^j \in \bigoplus_{k\in\wt'_\RC}\RC^k$ for $j > 0$ .
\end{proposition}
\pr
Since $R^0$ is the center of $R$, $L_{(k)}$ acts on $R/R^0$ for all $k\in \N$. So $R/R^0$ has the $sl_2$-module structure defined by
\begin{eqnarray}
 E &\mapsto& L_{(2)} ,\nonumber\\
 H &\mapsto& -2L_{(1)} ,\mylbl{decomp-1-1}\\
 F &\mapsto& -L_{(0)} .\nonumber\
\end{eqnarray}
Consider $P=\bigoplus_{k\in\wt'_R}{\Span\{x\in R/R^0|\ L_{(1)}x=kx,\ x\notin U,\ L_{(2)}x \in U\ \mbox{for some submodule}\ U\}}$.
Since $L_{(1)}P \subset P$, we have a basis $\{e_\lambda\}_{\lambda\in\Lambda}$ of $R/R^0$ and a function $\nu:\Lambda \to \wt'_R$ satisfying $L_{(1)}e_\lambda = \nu(\lambda) e_\lambda$ for all $\lambda \in \Lambda$ so that $\{e_\lambda\}_{\lambda\in\Lambda'}$ is a basis of $P$ for some $\Lambda'\subset \Lambda$.
Consider the $sl_2$-module homomorphism $ f: \bigoplus_{\lambda\in \Lambda'} V_\lambda \to R/R^0$ defined by $f(v_\lambda) = e_\lambda$ where $V_\lambda$ is the Verma module of the highest weight $-2\nu(\lambda)$ with respect to $H$ and $v_\lambda$ is the highest weight vector of $V_\lambda$.
Since $R$ is regular $f$ is surjective and
each Verma module $V_{-2\nu(\lambda)}$ is irreducible, hence $f$ is isomorphic.  Thus we have a unique decomposition $ x = \sum_j \partial^{(j)}x^j$ for any $ x\in R/R^0$ where $x^j \in P$ for all $j$.  Hence $P=\RC/R^0$, so we have the decomposition of the result.
The uniqueness is obvious.
\myqed

We shall call the decomposition of Proposition \ref{decomp-prop-1}  the {\it $\partial$-decomposition} of $x$ and $x^j$ the {\it $j$-part of $x$} setting  $x^j = 0$ for $\wt(x)-j \notin \wt_R$.

\addtocounter{mypcount}{1}\begin{corollary}\mylbl{decomp-cor-4}
$R=K[\partial]\RC$.
\end{corollary}

\addtocounter{mypcount}{1}\begin{corollary}\mylbl{decomp-cor-2}
$\RC$ is isomorphic to $R/(\partial R)$ as  $K$-vector spaces.
\end{corollary}

\addtocounter{mypcount}{1}\begin{corollary}
A regular conformal superalgebra $R$ over $\C$ is superconformal if and only if the reduced subspace is finite-dimensional, all the conformal weights are half-integers, the even subspace of the reduced subspace $\RC_{\mbox{even}} = \bigoplus_{n\in\N}\RC^n$ and the odd subspace $\RC_{\mbox{odd}} = \bigoplus_{n\in\N+\hf{1}}\RC^n$.
\end{corollary}

\addtocounter{mypcount}{1}\begin{corollary}\mylbl{decomp-cor-14}
 Let $\RC$ be the reduced subspace of a conformal superalgebra $R$ and  $ \RC_{(n)} $  a copy of $ \RC $ for each $n \in \Z$.
The Lie superalgebra associated to $ R $ is $ \bigoplus_{n\in \Z}{\RC_{(n)}} $ with the product
 $ [a_{(m)},b_{(n)}]=\sum_{j=0}^{\infty} {\left( {\matrix{m\cr
j\cr
}} \right)\left( {a_{(j)}b} \right)_{(n+m-j)}} $.
\end{corollary}

\addtocounter{mypcount}{1}\begin{proposition}\mylbl{decomp-prop-a}
For a homomorphism of conformal superalgebras\\ $f:(R,\{(n)\}_{n\in\N},L) \to (R',\{(n)\}_{n\in\N},L')$, {\rm (A)} $f(\RC) \subset \RC'$, {\rm (B)} $f(\RC) = \RC' $ if and only if $f$ is surjective, and {\rm (C)} $f|_{\RC}$ is injective if and only if $f$ is injective.
\end{proposition}
\pr 
{\rm (A)}  Since $0 = f(L_{(2)}a) = L'_{(2)}f(a)$ for $a\in\RC$, it is obvious $f(\RC) \subset \RC'$.

{\rm (B)}  Assume that $f(\RC) = \RC'$.  Consider the $\partial$-decomposition $ x' = \sum_j{\partial^{(j)}x'^j}$ for $x' \in R'$.  Then we have $ x^i \in \RC$ such that $f(x^i) = x'^i$ for all $i \in \N$, so $f(\sum_j{\partial^{(j)}x^i})=x$.  
Conversely assume that $f$ is surjective.  Then for any $ a' \in \RC'$ there exists $x \in R$ such that $f(x) = a'$, so $a' - \sum_j{\partial^{(j)} f(x^j)} = 0$.  Since the $\partial$-decomposition is unique, $f(x^0)=a'$ for some $x^0 \in \RC$.  Hence we have $f(\RC) = \RC'$.

{\rm (C)}  Assume that $f|_{\RC}$ is injective.  Take $x \in R$ such that $f(x) = 0$.  Then $ \sum_j{\partial^{(j)}f(x^j)} = 0$, so we have $f(x^i) = 0$ for all $i \in \N$.  Hence $x^i = 0$ for all $i\in \N$, which implies $x=0$.
The converse is obvious.
\myqed

\addtocounter{mypcount}{1}\begin{definition}{\rm \mylbl{decomp-df-21}
 Define the $\PP{n}$ product on $\RC$ for each $ n \in \N$ by
\begin{eqnarray*}
\PP{n}: \RC\times\RC &\to& \RC\\
 (a,b) &\mapsto& \pp{a}{n}{b} = (a_{(n)}b)^0,
\end{eqnarray*}
where $(a_{(n)}b)^0$ is the $0$-part of $a_{(n)}b$.
}\end{definition}

The center of a regular conformal superalgebra $R$ is $\{v\in\RC|\ \pp{v}{n}{x} = 0\ \mbox{for all}\ x \in \RC, n\in\N\}$.

\addtocounter{mypcount}{1}\begin{remark}{\rm 
The $\PP{n}$ products vanish except for finite many $\PP{n}$s if $\RC$ is finite-dimensional.
}\end{remark}

Let us denote $ (x;y)={{\Gamma (x+y)} \over {\Gamma (x)}}$ where $y$ is a non-negative integer and $x \in \C$, and define
\begin{eqnarray*}
G(\wt(a),\wt(b),n,j)&=&
\left\{
\begin{array}{c}
\begin{array}{r}
{{(2\wt(a)-n-j-1;j)} \over {(2(\wt(a)+\wt(b)-n-j-1);j)}}=\prod_{k=0}^{j-1}{\frac{(2\wt(a)-n-j-1+k)}{(2(\wt(a)+\wt(b)-n-j-1)+k)}} \\
\qquad \mbox{for}\ \wt(a)+\wt(b)-n-j-1 \notin -\hf{1}\N, 
\end{array}\\
\begin{array}{cl}
1, &\quad \mbox{for}\ \wt(a)+\wt(b)-n-1=0, j=0,\\
0, &\quad \mbox{otherwise}.
\end{array}
\end{array}
\right. \\
\end{eqnarray*}

\addtocounter{mypcount}{1}\begin{proposition}\mylbl{decomp-prop-22}
For a regular conformal superalgebra $R$
\begin{equation}
  (a_{(n)}b)^j=G(\wt(a),\wt(b),n,j)\pp{a}{n+j}{b}\mylbl{decomp-trans}
\end{equation}
where $a,b\in \RC$.
\end{proposition}
\pr
If $\wt(a)+\wt(b)-n-j-1=0$ the both sides are in $R^0$, so the proposition is obvious.  Otherwise,
apply $L_{(2)}^j$ to the both sides of $a_{(n)}b = \sum_j{\partial^{(j)}(a_{(n)}b)^j}$ and take the $0$-parts.   The left hand side is
\begin{equation}
  (L_{(2)}^ja_{(n)}b)^0=(-1)^j\left( {\prod\limits_{k=n}^{n+j-1} {\left( {k+2(1-\wt(a))} \right)}} \right)(a_{(n+j)}b)^0 .
\end{equation}
Taking the $0$-part of $ L^j_{(2)}\sum_j{\partial^{(j)}(a_{(n)}b)^j} $, we obtain
\begin{eqnarray}
{(L_{(2)}^ja_{(n)}b)^0}
 &=&\prod\limits_{k=1}^j {\left( {k+2(\wt(a)+\wt(b)-n-j-1)-1)} \right)}(a_{(n)}b)^j .
\end{eqnarray}
Hence we have
\begin{eqnarray}
  (a_{(n)}b)^j &=& \left( {\prod\limits_{k=0}^{j-1} { {{{-n-k+2\wt(a)-2} \over {k+2(\wt(a)+\wt(b)-n-j-1)}}} }} \right)\pp{a}{n+j}{b}\nonumber \\
  &=& G(\wt(a),\wt(b),n,j)\pp{a}{n+j}{b} .
\end{eqnarray}
\myqed

\addtocounter{mypcount}{1}\begin{proposition}\mylbl{decomp-prop-b}
A $K$-linear map $f: \RC \to \RC'$ satisfying $f(L) = L'$ and $f(\pp{a}{n}{b}) = \pp{f(a)}{n}{f(b)}$ for all $a,b \in \RC$, $n \in \N$ uniquely extends to a homomorphism of conformal superalgebras $\tilde{f}: R \to R'$, where $(R,\{(n)\}_{n\in\N},L)$ and $(R',\{(n)\}_{n\in\N},L')$ are regular conformal superalgebras.
\end{proposition}
\pr
Define a map $\tilde{f}: R \to R'$ by $\tilde{f}(x)=\sum_j{\partial^{(j)}f(x^j)}$ where $x\in R$.  Obviously $\tilde{f}(\partial x) = \partial \tilde{f}(x)$ for all $x \in R$.
By (\ref{decomp-trans}) we have
 $ {a_{(n)}b} = {\sum_{j=0}^\infty  {G(\wt(a),\wt(b),n,j)\partial ^{(j)}\pp{a}{n+j}{b}}} $,
 by (C1) and (C2)
$ (\partial^{(k)}a)_{(n)}(\partial^{(l)}b) = (-1)^k\sum_{j=0}^{l}{\frac{n!}{k! j! (n-k-j)!} \partial^{(l-j)}a_{(n-k-j)}b}$
for all  $a,b \in\RC$, $k,l\in\N$.  Hence the $(n)$ products on $R$ is written in terms of the $\PP{n}$ products and the operator $\partial$, so $\tilde{f}(x_{(n)}y)=\tilde{f}(x)_{(n)}\tilde{f}(y)$ for all $x,y\in R$ and $n\in\N$, that is, $\tilde{f}$ is an homomorphism of conformal superalgebras.

  Suppose given two extensions $\tilde{f}$ and $\tilde{f}'$ of $f$.  Then $\tilde{f}'(\partial^{(k)}a) = \partial^{(k)}f(a) = \tilde{f}(\partial^{(k)}a)$ holds for all $a\in \RC$.  Hence the extension of $f$ is unique by Corollary \ref{decomp-cor-4}.
\myqed 

\addtocounter{mypcount}{1}\begin{corollary}\mylbl{decomp-cor-c}
  Two conformal superalgebras $(R,\{(n)\}_{n\in\N},L)$ and $(R',\{(n)\}_{n\in\N},L')$ are isomorphic if and only if there exists a bijective $K$-linear map $f:\RC \to \RC'$ satisfying $f(L)=L'$ and $f(\pp{a}{n}{b}) = \pp{f(a)}{n}{f(b)}$ for all $a,b \in \RC$ and $n\in \N$.
\end{corollary}

We can reconstruct the ideals of a regular conformal superalgebra $(R,\{(n)\}_{n\in\N},L)$ from the triple $(\RC,\{\PP{n}\}_{n\in\N},L)$.

\addtocounter{mypcount}{1}\begin{proposition}\mylbl{primalg-prop-7}
 For an ideal $I$ of a conformal superalgebra $R$, there exists an ideal $\IC$ of the reduced subspace $\RC$ with respect to the $\PP{n}$ products.
 Conversely $I=K[\partial ]\IC$ is an ideal of $R$ for an ideal $\IC$ of $\RC$.
\end{proposition}
\pr
We may assume that $I$ is proper without loss of generality.
 Consider the projection $f: R \to R/I$.  We have a projection $\check{f}: \RC \to \check{(R/I)}$ with $\check{f}(\pp{a}{n}{b}) = \pp{\check{f}(a)}{n}{\check{f}(b)}$ for all $a,b \in \RC$. 
 Set $\IC = \ker \check{f}$.  We have $\pp{x}{n}{\IC}\subset\IC$ for all $x\in\RC$ and $I=\ker f = K[\partial]\IC$.  The converse is obvious.
\myqed

\addtocounter{mypcount}{1}\begin{corollary}\mylbl{primalg-cor-8}
A regular conformal superalgebra $(R,\{(n)\}_{n\in\N},L)$ is simple if and only if any ideal $\IC$ of the reduced subspace $\RC$ is either $\RC$ or $\{0\}$.
\end{corollary}

Consider the following properties of the triple $(\VC,\{\PP{n}\}_{n\in\N},L)$ for a $\ZZ$-graded $K$-vector space $\VC$ equipped with countably many products $\{\PP{n}\}_{n\in\N}$ on $V$ where $L\in\VC$:
\begin{namelist}{(P0) }
\item[(P0)] For $ a,b \in \VC $ there exists some  $ N \in \N $ such that for all $n \in \N$ satisfying $ n > N $,
 $$ \pp{a}{n}{b}=0 .$$
\item[(P2)] For $ a,b \in \VC $ and  $ n \in \N $,
 $$ \pp{a}{n}{b}=-(-1)^{n+p(a)p(b)}\pp{b}{n}{a} .$$ 
\item[(P3)] For $ a,b,c \in \VC $ and $ n,m \in \N $,
\begin{eqnarray*}
\lefteqn{  \sum^m\limits_{j=0} {\left( {\matrix{m\cr
j\cr
}} \right)G(\wt(b),\wt(c),n,j)\pp{a}{m-j}{\ppp{b}{n+j}{c}}}  } \\
& & -(-1)^{p(a)p(b)}\sum^n\limits_{j=0} {\left( {\matrix{n\cr
j\cr
}} \right)G(\wt(a),\wt(c),m,j)\pp{b}{n-j}{\ppp{a}{m+j}{c}}}\\
\qquad &=& \sum^{m+n}\limits_{j=0} {F(\wt(a),\wt(b),m,n,j) \pp{\left( \pp{a}{j}{b} \right)}{m+n-j}{c}}, 
\end{eqnarray*}
where
\begin{eqnarray*}
\lefteqn{ F(\wt(a),\wt(b),m,n,t)} \\
&=& \sum\limits_{k=0}^t {\left( {\matrix{m\cr
t-k\cr
}} \right)\left( {\matrix{{m+n+k-t}\cr
{k}\cr
}} \right)(-1)^{k}G(\wt(a),\wt(b),t-k,k)}.
\end{eqnarray*}
\item[(PV)] 
 $L$ is even and satisfies
 $\pp{L}{0}{a} = 0$,
 $\pp{L}{1}{L} = 2L$,
 $\pp{L}{2}{a} \in \VC^0$
for all  $ a\in \VC $.
The operator $ \pp{L}{1}{} $ is diagonalizable.  $\VC^0$ is central, $\wt_\VC \cap (-\hf{1}\N) \subset \{0\}$, and for all $k\in\wt_\VC$ there exists some $M\in\N$ such that $k-m \notin \wt_\VC$ for all $m\in\N$ satisfying $m\ge M$, where $\VC^k=\{a\in \VC|\ \pp{L}{1}{a}=ka\}$ and $\wt_\VC=\{k\in K|\ \VC^k\ne\{0\}\}$.
\end{namelist}

\addtocounter{mypcount}{1}\begin{proposition}\mylbl{decomp-prop-24}
The triple $(\VC,\{\PP{n}\}_{n\in\N},L)$ satisfies {\rm (P0)}, {\rm (P2)}, {\rm (P3)}, {\rm (PV)}, where
$\VC$ is the reduced subspace of a regular conformal superalgebra  with  the products $\{\PP{n}\}_{n\in \N}$.
\end{proposition}
\pr
Only (P3) is not obvious.
We shall obtain (P3) by taking the $0$-part of the both sides of (C3).
Apply Proposition \ref{decomp-prop-1} to the right hand side of (C3) and take the $0$-part.  Then,
\begin{eqnarray}
\lefteqn{ \left(\sum^m\limits_{k=0} {\left( {\matrix{m\cr
k\cr
}} \right)\left( {a_{(k)}b} \right)_{(m+n-k)}c}\right)^0  }\nonumber \\
 &=&\left(\sum^\infty\limits_{j,k=0} {\left( {\matrix{m\cr
 k\cr
 }} \right)\left( {\partial ^{(j)}(a_{(k)}b)^j} \right)_{(m+n-k)}c}\right)^0\nonumber \\
&=& \sum^\infty\limits_{j,k=0} {\left( {\matrix{m\cr
k\cr
}} \right)\left( {\matrix{{m+n-k}\cr
j\cr
}} \right)(-1)^jG(\wt(a),\wt(b),k,j)\pp{\left( {\pp{a}{k+j}{b}} \right)}{m+n-k-j}{c}}\nonumber \\
  &=&\sum^\infty\limits_{t=0} {\sum\limits_{k=0}^t {\left( {\matrix{m\cr
k\cr
}} \right)\left( {\matrix{{m+n-k}\cr
{t-k}\cr
}} \right)(-1)^{t-k}G(\wt(a),\wt(b),k,t-k)\pp{\left( {\pp{a}{t}{b}} \right)}{m+n-t}{c}}} \nonumber\\
 & =&\sum^{m+n}\limits_{t=0} {F(\wt(a),\wt(b),m,n,t) \pp{\left( \pp{a}{t}{b} \right)}{m+n-t}{c}},
\end{eqnarray}
which is the right hand side of (P3).
On the other hand the $0$-part of the term $a_{(m)}b_{(n)}c$ of the left hand side of (C3) is,
\begin{equation}
  \left(a_{(m)}b_{(n)}c\right)^0=\sum^m\limits_{j=0} {\left( {\matrix{m\cr
j\cr
}} \right)G(\wt(b),\wt(c),n,j)\pp{a}{m-j}{\ppp{b}{n+j}{c}}},
\end{equation} 
and for the term $b_{(n)}a_{(m)}c$
\begin{equation}
 \left(b_{(n)}a_{(m)}c\right)^0=\sum^n\limits_{j=0} {\left( {\matrix{n\cr
j\cr
}} \right)G(\wt(a),\wt(c),m,j)\pp{b}{n-j}{\ppp{a}{m+j}{c}}}.
\end{equation}
Thus we obtained the left hand side of (P3).
\myqed

\addtocounter{mypcount}{1}\begin{example}{\rm \mylbl{decomp-ex-27}
For $m=0$ and $n=0$, (P3) is
\begin{equation}
 \pp{a}{0}{\ppp{b}{0}{c}}-(-1)^{p(a)p(b)}\pp{b}{0}{\ppp{a}{0}{c}}=\pp{\left( {\pp{a}{0}{b}} \right)}{0}{c}. 
\end{equation}
For $m=1$ and $n=0$,
\begin{eqnarray}
 \lefteqn{\pp{a}{1}{\ppp{b}{0}{c}}+{{\wt(b)-1} \over {\wt(b)+\wt(c)-2}}\pp{a}{0}{\ppp{b}{1}{c}}-(-1)^{p(a)p(b)}\pp{b}{0}{\ppp{a}{1}{c}}}\nonumber \\
&=&\pp{\left( {\pp{a}{0}{b}} \right)}{1}{c}+{{\wt(b)-1} \over {\wt(a)+\wt(b)-2}}\pp{\left( {\pp{a}{1}{b}} \right)}{0}{c}.
\end{eqnarray}
For $m=1$ and $n=1$,
\begin{eqnarray}
\lefteqn{  \pp{a}{1}{\ppp{b}{1}{c}}+\frac{2\wt(b)-3}{2(\wt(b)+\wt(c)-3)}\pp{a}{0}{\ppp{b}{2}{c}}}\nonumber\\
\lefteqn{\qquad-(-1)^{p(a)p(b)}\left(\pp{b}{1}{\ppp{a}{1}{c}}+\frac{2\wt(a)-3}{2(\wt(a)+\wt(c)-3)}\pp{b}{0}{\ppp{a}{2}{c}}\right)}\nonumber\\
&=&\pp{\left(\pp{a}{0}{b}\right)}{2}{c} + {{\wt(b)-\wt(a)} \over {\wt(a)+\wt(b)-2}}\pp{\left( {\pp{a}{1}{b}} \right)}{1}{c}\nonumber \\
\lefteqn{- \frac{(2\wt(a)-3)(2\wt(b)-3)}{2(\wt(a)+\wt(b)-3)(2\wt(a)+2\wt(b)-5)}\pp{\left(\pp{a}{2}{b}\right)}{0}{c}.}
\end{eqnarray}
For $m=2$ and $n=0$,
\begin{eqnarray}
\lefteqn{ \pp{a}{2}{\ppp{b}{0}{c}} +   {{2(\wt(b)-1)} \over {\wt(b)+\wt(c)-2}}\pp{a}{1}{\ppp{b}{1}{c}}}\nonumber\\
\lefteqn{\qquad + \frac{(2\wt(b)-3)(\wt(b)-1)}{(\wt(b)+\wt(c)-3)(2\wt(b)+2\wt(c)-5)}\pp{a}{0}{\ppp{b}{2}{c}}}\nonumber\\
\lefteqn{\qquad - (-1)^{p(a)p(b)}\pp{b}{0}{\ppp{a}{2}{c}}}\nonumber \\
&=& \pp{\left(\pp{a}{0}{b}\right)}{2}{c} { + {{2(\wt(b)-1)} \over {\wt(a)+\wt(b)-2}}\pp{\left( {\pp{a}{1}{b}} \right)}{1}{c}}\nonumber\\
\lefteqn{  +\frac{(\wt(b)-1)(2\wt(b)-3)}{(\wt(a)+\wt(b)-3)(2\wt(a)+2\wt(b)-5)}\pp{\left(\pp{a}{2}{b}\right)}{0}{c}.}
\end{eqnarray}
}\end{example}

\section{ Reconstruction of the conformal superalgebras }
We can reconstruct the entire regular conformal superalgebra $(R,\{(n)\}_{n\in\N},L)$ from the triple $(\RC,\{\PP{n}\}_{n\in\N},L)$.

\addtocounter{mypcount}{1}\begin{theorem}\mylbl{primalg-th-3}
For a triple $(\VC,\{\PP{n}\}_{n\in\N},L)$ satisfying {\rm (P0)}, {\rm (P2)}, {\rm (P3)} and {\rm (PV)}, there exists a regular conformal superalgebra $ (R_\VC,\{(n)\}_{n\in\N},L) $ whose reduced subspace is $\VC$ and the products satisfies $(a_{(n)}b)^0 = \pp{a}{n}{b}$ for all $a,b\in \VC$, $n\in\N$.  Furthermore the conformal superalgebra is unique up to isomorphisms.
\end{theorem}
\pr
Consider the left $K[\partial]$-module 
$R_\VC = \left(K[\partial] \otimes (\bigoplus_{k\in(\wt_\VC\setminus\{0\})} \VC^k)\right) \oplus \VC^0, $
 where $\partial$ is an indeterminate and $\VC^0$ is regarded as a left $K[\partial]$-module by $\partial \VC^0 =0$.  We  omit the $\otimes$ for brevity.
Define $ \wt(\partial^j a)=\wt(a)+j $ for $a \in \VC$.
 Each  $ x\in R_\VC $ is uniquely  written as
 $ x = \sum_{j}{\partial^{(j)}x^j} $
 for some $x^0 \in \VC$ and $ x^j \in \bigoplus_{k\in(\wt_\VC\setminus\{0\}}\VC^k $ for all $j>0$.
 Define  $(n)$ products on  $ R_\VC $ by
 $ a_{(n)}b=\sum_{j=0}^\infty  {G(\wt(a),\wt(b),n,j)\partial ^{(j)}\pp{a}{n+j}{b}} $
and
$ (\partial^{(k)}a)_{(n)}(\partial^{(l)}b) = (-1)^k\sum_{j=0}^{l}{\frac{n!}{k! j! (n-k-j)!} \partial^{(l-j)}a_{(n-k-j)}b}$ where $ a,b \in \VC $.
It is easy to check that the $(n)$ products satisfy (C0) and  (C1), 
and by direct calculation $\partial$ is a derivation with respect to the $(n)$ products.
Now, for $a,b \in \VC$,
\begin{eqnarray}
\lefteqn{ -(-1)^{p(a)p(b)+n} \sum_{j=0}^\infty{(-1)^j\partial^{(j)}(b_{(n+j)}a)}}\nonumber\\
&=& -(-1)^{p(a)p(b)+n} \sum_{j,k=0}^\infty{(-1)^j\partial^{(j)}G(\wt(b),\wt(a),n+j,k)\partial^{(k)}(\pp{b}{n+j+k}{a})}\nonumber\\
&=& \sum_{s=0}^\infty \sum_{k=0}^{s}{(-1)^{s-k+p(a)p(b)+n+1}\left(\matrix{s\cr k}\right) G(\wt(b),\wt(a),n+s-k,k)\partial^{(s)}(\pp{b}{n+s}{a})}\nonumber\\
&=&  \pp{a}{n}{b} + \sum_{s=1}^\infty {}_{2}F_{1}\left(\matrix{\matrix{-s & 2\wt(b)-n-s-1}\cr 2(\wt(a)+\wt(b)-n-s-1)};1\right)\partial^{(s)}(\pp{a}{n+s}{b}),
\end{eqnarray}
where $_{2}F_{1}\left(\matrix{\matrix{\alpha& \beta} \cr \gamma }; x \right) = \sum_{j=0}^\infty {\frac{(\alpha;j)(\beta;j)}{(\gamma;j)}x^{(j)}}$.
We have $_{2}F_{1}\left(\matrix{\matrix{\alpha & \beta}\cr \gamma} ; 1 \right) = \frac{\Gamma(\gamma)\Gamma(\gamma-\alpha-\beta)}{\Gamma(\gamma-\alpha)\Gamma(\gamma-\beta)}$ for $\gamma\notin -\N$ and $\alpha\in -\N$.  If $s\ge 1$ and $\wt(a)+\wt(b)-n-s-1 = \wt(\pp{a}{n+s}{b}) \in -\hf{1}\N$ then $\partial^{(s)}(\pp{a}{n+s}{b}) = 0$, so we have
\begin{eqnarray}
\lefteqn{ -(-1)^{p(a)p(b)+n} \sum_{j=0}^\infty{(-1)^j\partial^{(j)}(b_{(n+j)}a)}=\pp{a}{n}{b}}& &\nonumber\\
& & +\sum_{s=1}^\infty \frac{\Gamma(2(\wt(a)+\wt(b)-n-s-1))\Gamma(2\wt(a)-n-1)}{\Gamma(2(\wt(a)+\wt(b)-n-s-1)+s)\Gamma(2\wt(a)-n-s-1)}\partial^{(s)}(\pp{a}{n+s}{b})\nonumber\\
&=&  \sum_{s=0}^\infty G(\wt(a),\wt(b),n,s)\partial^{(s)}(\pp{a}{n+s}{b})\nonumber\\
&=&  a_{(n)}b.
\end{eqnarray}
Then, (C2) is checked for all $a=\partial^kx$ and $b=\partial^ly$ by induction on $k$ and $l$ where $x,y \in \VC$.
Indeed, assume $(\partial^kx)_{(n)}(\partial^ly) = (-1)^{p(x)p(y)}\sum_j (-1)^{1+n+j}\partial^{(j)}(\partial^ly)_{(n+j)}(\partial^kx)$ for all $n$.  Applying $\partial$ to the both sides we have
$(\partial^{k+1}x)_{(n)}(\partial^ly) + (\partial^kx)_{(n)}(\partial^{l+1}y)
 = (-1)^{p(x)p(y)}\sum_j{\frac{(-1)^{1+n+j}}{j!}\partial^{j+1}(\partial^ly)_{(n+j)}(\partial^kx)}$.
By (C1) we have
\begin{eqnarray}
\lefteqn{(\partial^kx)_{(n)}(\partial^{l+1}y)}\nonumber\\
 &=& n(\partial^kx)_{(n-1)}(\partial^ly) + (-1)^{p(x)p(y)}\sum_j{(-1)^{n+j}j\partial^{(j)}(\partial^ly)_{(n+j-1)}(\partial^kx) }\nonumber\\
&=& (-1)^{p(x)p(y)}\sum_j{(-1)^{1+n+j}\partial^{(j)}(\partial^{l+1}y)_{(n+j)}(\partial^kx)},
\end{eqnarray}
 which implies (C2) for $a=\partial^kx$ and $b=\partial^{l+1}y$. 
On the other hand
\begin{eqnarray}
& &\!\!\!\!\!\!\!\!(\partial^{k+1}x)_{(n)}(\partial^ly) + (\partial^kx)_{(n)}(\partial^{l+1}y) \nonumber\\
& &\!\!\!\!\!\!\!\!\!\!\!\!\!\!\!\!\! =(-1)^{p(x)p(y)}\sum_j{ (-1)^{1+n+j}\partial^{(j)}\left((\partial^{l+1}y)_{(n+j)}(\partial^kx) + (\partial^ly)_{(n+j)}(\partial^{k+1}x) \right) },
\end{eqnarray}
 thus we have (C2) for $a=\partial^{k+1}x$ and $b=\partial^ly$, which completes the induction.

For all $a \in \VC$, we have
\begin{eqnarray}
 L_{(0)}(\partial^k a) &=& \partial^{k+1} a ,\mylbl{prim-3-1}\\
 L_{(1)}(\partial^k a) &=& (\wt(a)+k) \partial^k a \mylbl{prim-3-2},\\
 L_{(2)}(\partial^k a) &=& k(k-1+2\wt(a)) \partial^{k-1} a ,\mylbl{prim-3-3}
\end{eqnarray}
which imply (CV).

In order to show (C3), let
\begin{eqnarray}
\lefteqn{\!\!\!J(a,b,c,m,n,k)}\nonumber\\
& &\!\!\!\!\!\!\!\! = a_{(m)}b_{(n)}\partial^{k}c - (-1)^{p(a)p(b)}b_{(n)}a_{(m)}\partial^{k}c - \sum_{j=0}^{\infty}{\left( {\matrix{m\cr j}}\right)(a_{(j)}b)_{(m+n-j)}\partial^{k}c},
\end{eqnarray}
where $a,b,c\in\VC$.  By (C1) and the Leibniz rule we have 
\begin{eqnarray}
\lefteqn{\partial J(a,b,c,m,n,k) }\nonumber\\
& & = -m J(a,b,c,m-1,n,k) -n J(a,b,c,m,n-1,k)  + J(a,b,c,m,n,k+1),\mylbl{prim-3-4}
\end{eqnarray}
 where we understand $J(a,b,c,-1,n,k)=0$ and $J(a,b,c,m,-1,k)=0$ for all $m,n,k \in \N$.
We have  $\wt(b_{(n)}\partial^kc) - \wt(\partial^kc) = \wt(b)-n-1$, so  $J(L,b,c,1,n,k)=0$ for all $n,k\in \N$.  By (P3) and the definition of $L_{(1)}$, $J(L,b,c,2,n,0) = 0$.  Substituting them into (\ref{prim-3-4}) for $m=2$ we obtain $J(L,b,c,2,n,k)=0$ for all $n,k \in \N$ by induction on $k$.
On the other hand taking the $0$-part of (\ref{prim-3-4}) we obtain $0 = - m J(a,b,c,m-1,n,k)^0 - n J(a,b,c,m,n-1,k)^0 + J(a,b,c,m,n,k+1)^0$.  By induction on $k$ we have  $J(a,b,c,m,n,k)^0 =0$  for any $m,n,k \in \N$.

Consider
\begin{equation}
  B(a,b,c,m,n)= {a_{(m)}b_{(n)}c-(-1)^{p(a)p(b)}b_{(n)}a_{(m)}c-\sum\limits_{j=0}^\infty  {\left( \matrix{m\hfill\cr
j } \right)\left( {a_{(j)}b} \right)_{(m+n-j)}c}},
\end{equation}
where $a,b,c\in R_\VC$.
By (C1) and (C2), 
it suffices to check $B(a,b,c,m,n)=0$ for $m,n \in \N$ where $a,b \in \VC$ and $c \in R_\VC$.  We have $B(a,b,c,m,n)^0=0$ for $a,b \in \VC$ and $c\in R_\VC$ because $J(a,b,c,m,n,k)^0 =0$ for all $k\in\N$.
If $\wt(B(a,b,c,m,n))=0$ then $B(a,b,c,m,n)=B(a,b,c,m,n)^0=0$,
so we may assume $\wt(B(a,b,c,m,n))\ne 0$.
By (\ref{prim-3-3}) we have
\begin{equation}
  \left((L_{(2)})^kB(a,b,c,m,n)\right)^0=(2(\wt(a)+\wt(b)+\wt(c)-m-n-2-k);k)\left(B(a,b,c,m,n)\right)^k.
\end{equation}
The coefficients on the right hand side never vanish because $ \wt(B(a,b,c,m,n)^k)\notin -\hf{1}\N $, hence $B^k$ is proportional to $((L_{(2)})^k B)^0$.
 Since $J(L,b,c,2,n,k)=0$,
\begin{eqnarray}
\lefteqn{L_{(2)}B(a,b,c,m,n)}\nonumber\\ &=& -(n+2(1-\wt(a)))B(a,b,c,m+1,n) \nonumber\\
& &\  -(m+2(1-\wt(b)))B(a,b,c,m,n+1) + B(a,b,L_{(2)}c,m,n),\mylbl{prim-3-5}
\end{eqnarray}
$ (L_{(2)})^k B(a,b,c,m,n) $ is written by a linear combination of some $ B $s.
  Thus we have $ B^k=0 $ for all $ k\in \N $, which implies (C3).

The reduced subspace of $R_\VC$ coincides with $\VC$ itself.
The uniqueness follows from Corollary \ref{decomp-cor-c}.
\myqed

The following lemma plays an important role in later sections.

\addtocounter{mypcount}{1}\begin{lemma}\mylbl{primalg-lem-solve}
\begin{equation}
  \pp{a}{p}{\ppp{b}{q}{c}}=\sum\limits_j {r_{j}\pp{b}{q-j}{\ppp{a}{p+j}{c}}+s_{j}\pp{\left( {\pp{a}{j}{b}} \right)}{p+q-j}{c}} \mylbl{eqsolve}
\end{equation}
for some $r_{j},s_{j} \in K$ where $a,b,c \in \RC$, $ p,q \in \N$.
\end{lemma}
\pr
Denote $\mbox{(P3)}{}_{m,n}$ for (P3) specifying $m,n$.
For $p=0$, (\ref{eqsolve}) follows from $\mbox{(P3)}{}_{0,q}$.  
Suppose (\ref{eqsolve}) holds for all $p\le k$ and $q\in\N$.  $\mbox{(P3)}{}_{k+1,q}$ implies (\ref{eqsolve}) for $p=k+1$.
\myqed

\addtocounter{mypcount}{1}\begin{proposition}\mylbl{primalg-prop-9}
Consider a regular conformal superalgebra
 $(R,\{(n)\}_{n\in\N},L)$.
Let $S$ be a subset of the reduced subspace $\RC$ and $I_S$ the ideal generated by $S$.
Then for a basis $B$ with an order $<$ on $B$,
\begin{equation}
 I_S\cap \RC = \Span\left\{\pp{v^1}{n_1}{\ppp{v^2}{n_2}{\cdots \ppp{v^r}{n_r}{u}}}\left|\ {v^k\in B,\ u\in S,\ n_k\in \N,\ v^i< v^{i+1}} \right.\right\}.\label{ordered}
\end{equation}
\end{proposition}
\pr
Set
\begin{equation}
 F_p\IC_S = \Span\left\{\pp{v^1}{n_1}{\ppp{v^2}{n_2}{\ppp{v^3}{n_3}{\cdots \ppp{v^r}{n_r}{u}}}}\left| {v^k,u\in S,\ n_k \in \N,\ r \leq p} \right.\right\}.
\end{equation}
By Lemma \ref{primalg-lem-solve},
 $ \pp{v^1}{n_1}{\ppp{v^2}{n_2}{\ppp{v^3}{n_3}{\cdots\ppp{v^r}{n_r}{u}}}} $ 
is written by a linear combination of these elements with $v^i$ and $v^{i+1}$ swapped (but $n_l$s may differ) as an element of $ (F_{i}\IC_S) / (F_{i-1}\IC_S)  $.
Hence we have
\begin{equation}
 F_p\IC_S = \Span\left\{\pp{v^1}{n_1}{\ppp{v^2}{n_2}{\cdots \ppp{v^r}{n_r}{u}}}\left|\ {v^k\in B,\ u\in S,\ n_k\in \N,\ r\leq p,\ v^i< v^{i+1}} \right.\right\}.
\end{equation}
Since $I_S\cap \RC = \sum_p F_p\IC_S$, thus we have the result.
\myqed

Let $\CP$ be the category of triples $(\VC,\{\PP{n}\},L)$ satisfying (P0), (P2), (P3), (PV) where $\VC$ is a vector space, $\{\PP{n}\}$ is a set of products on $\VC$, and $L$ is a vector in $\VC$, with the morphisms being the linear maps that commute with all the $\PP{n}$ products and preserve $L$.  We can summarize this section:  the category of regular conformal superalgebras is equivalent to the category $\CP$ by the functor $F(R) = \RC$ and $F(f) = f|_{\RC}$.

\section{ Physical conformal superalgebra}

In this section we will study physical conformal superalgebras.  One can  reduce the axioms for conformal superalgebras into some simple relations.   We shall assume $K=\C$ hereafter.

A regular conformal superalgebra $R$ is physical if and only if the reduced subspace $\RC$ satisfies the following.
\begin{itemize}
\item Eigenvalues of $ L_{(1)} $ on $\RC$ are $2$, $\hf{3}$, $1$ and $\hf{1}$.
\item $\RC^2 = \K L $.
\item $\RC^{3/2}  $ and $\RC^{1/2}$ are odd subspaces.
\item $\RC^{1} $ and $\RC^{2}$ are even subspaces.
\end{itemize}
All the $\PP{n}$ products vanish except for the $\PP{0}$ product and the $\PP{1}$ product for physical conformal superalgebras.

Let $R$ be a physical conformal superalgebra and $\RC$ be the reduced subspace. Consider the products on $\RC$ defined by,
 $$  a \mdot b=\left\{ 
\begin{array}{cc}
 \frac{\pp{a}{1}{b}}{\wt(a)+\wt(b)-2}, &\quad \mbox{for } \wt(a)+\wt(b)-2 \ne 0 ,\\
 0 ,& \quad \mbox{otherwise}, 
\end{array}
\right. $$ 
 $$ a\bdot b=\pp{a}{0}{b}. $$ 
Obviously we have
$$ L \mdot a = a,\  L\bdot a =0,  \leqno{\mbox{(D0)}}$$
 $$  a\mdot b = (-1)^{p(a)p(b)} b\mdot a , \leqno{\mbox{(D1)}}$$ 
 $$  a\bdot b = -(-1)^{p(a)p(b)} b\bdot a . \leqno{\mbox{(D2)}} $$ 
Rewriting the relations in Example \ref{decomp-ex-27} in terms of the product $\mdot$ and the product $\bdot$ we obtain
 $$  (\wt(b)-1)a\mdot b\mdot c = (\wt(b)-1)(a \mdot b) \mdot c , \leqno{\mbox{(D3)}}  $$ 
 $$  (\wt(b)+\wt(c)-2)a\mdot b\mdot c-(-1)^{p(a)p(b)}(\wt(a)+\wt(c)-2)b\mdot a\mdot c  \leqno{\mbox{(D4)}}$$
$$=(\wt(b)-\wt(a))\left( {a\mdot b} \right)\mdot c , $$ 
 $$  \left( {\wt(a)+\wt(b)+\wt(c)-3} \right)a\mdot b\bdot c+\left( {\wt(b)-1} \right)a\bdot b\mdot c \leqno{\mbox{(D5)}}$$
$$-(-1)^{p(a)p(b)}\left( {\wt(a)+\wt(c)-2} \right)b\bdot a\mdot c  $$ 
 $$ =\left( {\wt(a)+\wt(b)+\wt(c)-3} \right)\left( {a\bdot b} \right)\mdot c+\left( {\wt(b)-1} \right)\left( {a\mdot b} \right)\bdot c ,  $$ 
 $$  a\bdot b\bdot c-(-1)^{p(a)p(b)}b\bdot a\bdot c=\left( {a\bdot b} \right)\bdot c .\leqno{\mbox{(D6)}}$$

Let us denote  $ V =\RC^{3/2}$, $A= \RC^1$, $ F= \RC^{1/2}$ for a physical conformal superalgebra $R$, following the notations in \cite{kac2}.
That is,  $\RC$ is decomposed into
$ \RC= \K L \oplus V \oplus A \oplus F$.
Define the inner product $(\cdot,\cdot)$ on $V$ by $\pp{a}{0}{b}=(a,b)L$.
Consider the following properties:
\begin{namelist}{(Hx) }
\item[(H0)] ${L\mdot}= \mbox{\it id}$, ${L\bdot} = 0$ as operators on $\RC$,
\item[(H1)] $(\RC,\bdot)$ is a Lie superalgebra,
\item[(H2)] $(\RC,\mdot)$ is an associative commutative superalgebra,
\item[(H3)] $A\bdot$ gives derivations with respect to $\mdot$,
\item[(H4)] $ u \mdot v \bdot f = (u\mdot v)\bdot f + (u \bdot v)\mdot f$, for $u,v \in V$ and $f \in F$,
\item[(H5)] $ ( u \mdot + u \bdot )^2 v =  (u,u) v,$ for $u,v \in V$,
\item[(H6)] $  ( u \mdot + u \bdot )^2 a =  (u,u) a,$ for $ u \in V$ and $ a \in A$.
\end{namelist}

\addtocounter{mypcount}{1}\begin{proposition}\mylbl{dotalg-prop-8.1}
 For the reduced subspace $\RC = \K L \oplus V \oplus A \oplus F$  of a physical conformal superalgebra, the products $\mdot$ and $\bdot$ have the properties {\rm (H0-6)}.
\end{proposition}
\pr
(D2) and (D6) yield (H1).
The product $\mdot$ is commutative by (D1).
We have
$ (\wt(b)-1)a\mdot b \mdot c = (-1)^{p(a)p(c)+p(b)p(c)}(\wt(b)-1)c\mdot a\mdot b$ by (D1) and (D3).
If $\wt(a)=\wt(b)=\wt(c)=1$ then the both sides are $0$ because $\wt(a)+\wt(b)+\wt(c)-4 = -1$.  So we may assume $\wt(b) \ne 1$ without loss of generality, hence $a\mdot b \mdot c = (-1)^{p(a)p(c)+p(b)p(c)}c\mdot a\mdot b$, thus $a\mdot b \mdot c = (a\mdot b)\mdot c$, we have (H2).
For the others, let
\begin{eqnarray}
Q(x,y,z) &=& (\wt(x)+\wt(y)+\wt(z)-3)x\mdot y\bdot z +(\wt(y)-1)x\bdot y\mdot z \nonumber\\
& &  -(-1)^{p(x)p(y)}(\wt(x)+\wt(z)-2)y\bdot x\mdot z - (\wt(y)-1)(x\mdot y)\bdot z \nonumber\\
& & - (\wt(x)+\wt(y)+\wt(z)-3)(x\bdot y)\mdot z,\mylbl{transQ}
\end{eqnarray}
and
\begin{eqnarray}
P(x,y,z) &=& (\wt(x)-1)x \mdot y \bdot z + (\wt(y) -1)x \bdot y \mdot z - (-1)^{p(x)p(y)}(\wt(y) -1) y \mdot x \bdot z\nonumber\\
 & & - (-1)^{p(x)p(y)}(\wt(x)-1)y \bdot x \mdot z - (\wt(x)+\wt(y) -2)(x \bdot y) \mdot z.\mylbl{transP}
\end{eqnarray}
It is easy to check 
$(\wt(x)-1)Q(x,y,z) + (\wt(y)-1)(-1)^{p(x)p(y)}Q(y,x,z) = (\wt(x)+\wt(y)+\wt(z)-3)P(x,y,z)$.
Since if $\wt(x)+\wt(y)+\wt(z)-3 = 0$ then the both sides are $0$, so we have $P(x,y,z)=0$ for all $x,y,z \in \RC$ because $Q(x,y,z)=0$ for all $x,y,z\in\RC$ by (D5).
$P(x,y,z)=0$ for $\wt(z)=1$, $(\wt(x),\wt(y),\wt(z)) = (\hf{1},\hf{3},\hf{3})$, $(\hf{3},\hf{3},\hf{3})$ and $(\hf{3},\hf{3},1)$ imply (H3), (H4), (H5) and (H6) respectively.\\
\myqed

 (H4), (H5) and (H6) imply the following.

\addtocounter{mypcount}{1}\begin{proposition}\mylbl{dotalg-prop-10}
The reduced subspace $\RC$ of a physical conformal superalgebra is a left $\Cl(V,(\cdot,\cdot))$-module by the action
$vx = v \mdot x + v \bdot x,$
where $ v \in V$ and $x \in \RC$.
\end{proposition}

Thus we have obtained the action of the Clifford algebra $\Cl(V,(\cdot,\cdot))$ on the associated Lie superalgebra, where $V$ is the space of the reduced vectors with the conformal weight $\hf{3}$.   The action is discussed in \cite{HTT}.

\addtocounter{mypcount}{1}\begin{corollary}\mylbl{dotalg-cor-11}
The Clifford algebra $\Cl(V,(\cdot,\cdot))$ acts on the associated Lie superalgebra of a physical conformal superalgebra, where $V$ is the space of reduced vectors of the conformal weight $\hf{3}$ with the inner product defined by $(u,v)L = u_{(0)}v$.
\end{corollary}

Furthermore we have the converse of Proposition \ref{dotalg-prop-8.1}.
\addtocounter{mypcount}{1}\begin{proposition}\mylbl{dotalg-prop-8.1a}
 Suppose given a finite-dimensional $\ZZ$-graded vector space $\RC$ 
with the decomposition $\RC = \K L \oplus V \oplus A \oplus F$ with respect to a weight $\wt$, where $\wt(V)=3/2$, $\wt(A)=1$, and $\wt(F)=1/2$ with the parity $p(\K L)= p(A)=0$ and $p(V)=p(F)=1$,
and two products $\mdot$ and $\bdot$ with the weight
$\wt(x\bdot y) = \wt(x) + \wt(y) -1$, $\wt(x\mdot y) = \wt(x) + \wt(y) -2$.
If $(\RC,\bdot,\mdot)$ have the properties {\rm (H0-6)} then 
the triple $(\RC,\{\PP{n}\},L)$ is a physical conformal superalgebra where we set $\pp{a}{0}{b} = a\bdot b$, $\pp{a}{1}{b} = (\wt(a)+\wt(b)-2)a\mdot b$ and $\pp{a}{n}{b}=0$ for $n\ge 2$.
\end{proposition}
\pr
  (P0), (P2) and (PV) are obvious.  It is easy to check that $\mbox{(P3)}{}^{a,b}_{m,n}$ is equivalent to $\mbox{(P3)}{}^{b,a}_{n,m}$, hence {\rm (D0-6)} are sufficient to (P3).  Only (D5) is not obvious since (D0) is (H0) itself, (D1) and (D6) follow from (H1), and  (D2), (D3) and  (D4) follow from (H2).
Let $P$ and $Q$ be as in (\ref{transQ}) and (\ref{transP}).
Then, (H1) and (H2) imply
\begin{equation}
 -(-1)^{p(y)p(z)}P(x,z,y) + (-1)^{p(x)p(y)+p(x)p(z)}P(y,z,x) =(-1)^{p(y)p(z)} Q(x,z,y),
\end{equation}
hence if $P(x,y,z)=0$ for all $x,y,z \in \RC$ then we have $Q(x,y,z)=0$ for all $x,y,z \in \RC$, which implies (D5).
Let us show $P(x,y,z)=0$ for all $x,y,z\in\RC$.
Since $P$ satisfies $P(x,y,z) = -(-1)^{p(x)p(y)}P(y,x,z)$ and 
${(-1)^{p(x)p(y)}(\wt(z)-1)P(x,y,z) - (-1)^{p(x)p(y)+p(y)p(z)}(\wt(y)-1)P(x,z,y)}$\\
$ = - (-1)^{p(x)p(z)}(\wt(x)-1) P(y,z,x)$,
 we have
\mbox{(1)} $P(x,y,z)=0 \Leftrightarrow P(y,x,z)=0$ and
\mbox{(2)} if $\wt(x) \ne 1 $, $ P(x,y,z)=0 \wedge P(x,z,y)=0 \Rightarrow P(y,z,x)=0$.
So it suffices to consider the following cases: $(\wt(x),\wt(y),\wt(z)) = (\hf{3},\hf{3},\hf{3})$,
$(\hf{3},1,\hf{3})$,
$(\hf{3},1,1)$,
$(\hf{3},\hf{3},1)$,
$(\hf{1},\hf{3},\hf{3})$.
 (H3) gives $(\hf{3},1,\hf{3})$ and $(\hf{3},1,1)$.
(H5), (H6) and (H4) imply $(\hf{3},\hf{3},\hf{3})$, $(\hf{3},\hf{3},1)$ and $(\hf{1},\hf{3},\hf{3})$ respectively.  Hence (D5) is shown for all $x,y,z \in \RC$.
\myqed

\section{ Simple physical conformal superalgebra}

In this section we will describe some properties of simple physical conformal superalgebra.  A criterion for simplicity is given.
Let $R$ be a simple physical conformal superalgebra, $\RC$ the reduced subspace, $V=\RC^\hf{3}$, $A=\RC^1$ and $F=\RC^\hf{1}$.

The following result is stated in \cite{kac1} and \cite{kac2}:  
\addtocounter{mypcount}{1}\begin{proposition}\mylbl{simplephy-prop-1}
Let $\RC$ be the reduced subspace of a simple physical conformal superalgebra. Then the inner product $(\cdot,\cdot)$ on $V$ is nondegenerate. 
\end{proposition}
\pr
Set $ V^0=\left\{ {v\in V\left| {\mbox{\ for all\,} u\in V \ (u,v)=0} \right.} \right\}$ and consider the ideal $I_{V^0}$ generated by $V^0$.
Fix a basis $ B $ of $\RC$ and take  an order on $B$ such that  $ a <b $  if $ \wt(a) < \wt(b) $.  Set\\
 $ F_p\JC=\Span\left\{{ \pp{x^1}{i_1}{\ppp{x^2}{i_2}{\cdots \ppp{x^k}{i_k}{v_0}}}\left| {x^r\in B,\ v_0\in V^0,\ x^i \leq x^{i+1},\ k\leq p} \right.} \right\}$.  Obviously $L\notin F_0\JC$.
 Take any $ u,v \in V $, $ v_0 \in V^0 $, $ x^r \in B $.
Since $ v\bdot v_0=0$, we have
 $ \pp{x^1}{i_1}{\cdots \ppp{x^{p-1}}{i_{p-1}}{ v\bdot v_0}}= 0 $ 
in $ F_p\JC/F_{p-1}\JC $.
By (D5) we have
\begin{eqnarray}
 \pp{x^1}{i_1}{\cdots \ppp{x^{p-2}}{i_{p-2}}{ v\bdot u\mdot v_0}}&=& \pp{x^1}{i_1}{\cdots \ppp{x^{p-2}}{i_{p-2}}{ (- \hf{1} u\bdot v \mdot v_0 )}}\nonumber \\
  &=& \pp{x^1}{i_1}{\cdots \ppp{x^{p-2}}{i_{p-2}}{ (\frac{1}{4} v \bdot u \mdot v_0)}}\nonumber\\
 &=& 0,
\end{eqnarray}
in $ F_p\JC/F_{p-1}\JC $.
By (D3) we have
 $\pp{x^1}{i_1}{\cdots \ppp{x^{p-2}}{i_{p-2}}{ v\mdot u\mdot v_0}}= 0$
in $ F_p\JC/F_{p-1}\JC $.
 $ \wt(\pp{x}{n}{y}) > \wt(y) $ occurs  only when $ x \in V$, thus we have $L\notin \sum_p{F_p\JC}$.
By Proposition \ref{primalg-prop-9} $\sum_p{F_p\JC} = I_{V^0}\cap \RC$, hence $I_{V^0}$ is proper unless $V^0 = \{0\}$.
\myqed

Now, let
$ F^3 = \{ f \in F | \pp{v^1}{0}{\ppp{v^2}{0}{\ppp{v^3}{0}{f}}} = 0\ \mbox{for all}\  v^k \in V \}$.

\addtocounter{mypcount}{1}\begin{proposition}\mylbl{simplephy-prop-4}
A physical conformal superalgebra $R$ with $V\ne \{0\}$ is simple if and only if $F^3 = 0$ and the inner product is nondegenerate.
\end{proposition}
\pr
Suppose that the inner product on $V\ne\{0\}$ is nondegenerate and $R$ is not simple. Take $\IC$ a proper ideal of the reduced subspace $\RC$.  $\IC$ is decomposed into $\IC = (\IC\cap \K L) \oplus (\IC\cap V) \oplus (\IC\cap A) \oplus(\IC\cap F)$ by the action of $\pp{L}{1}{}$.  $\IC\cap \K L = \{0\}$ because $\IC \ne \RC$.  In particular $\IC\cap V = \{0\}$ since the inner product on $V$ is nondegenerate.  $\IC$ is a $\Cl(V,(\cdot,\cdot))$-module because $I$ is an ideal.
The Clifford action of a unit vector in $V$ yields an isomorphism of vector spaces between $\IC\cap (\K L \oplus A)$ and $ \IC\cap (V\oplus F)$, so if $\IC \cap F = \{0\}$ then $\IC=\{0\}$.  Since $L\notin \IC$ we have $\IC \cap F \subset F^3$ while $I \ne \{0\}$, thus $F^3 \ne \{0\}$.

Conversely assume $F^3 \ne \{0\} $.  Let $I$ be the ideal generated by $F^3$. 
Apply Proposition \ref{primalg-prop-9} for $S=F^3$ taking an order such that $x < y $ if $\wt(x) < \wt(y)$. 
Take $f \in F^3$, $ v^i \in V ,(i = 1,2,3)$, and $ a \in A$.  Then,
\begin{eqnarray}
 v^1 \bdot v^2 \bdot v^3 \bdot a \bdot f &=& a \bdot v^1 \bdot v^2 \bdot v^3 \bdot f - (a\bdot v^1) \bdot v^2 \bdot v^3 \bdot f\nonumber\\
& &  - v^1 \bdot (a\bdot v^2)\bdot v^3 \bdot f - v^1\bdot v^2\bdot (a \bdot v^3)\bdot f\nonumber\\
&=& 0,
\end{eqnarray}
so $ A \bdot F^3 \subset F^3$.  
Applying (D5) for any $a,b \in V $ and $ c \in F^3 $, we have 
$a \mdot b \bdot c = (a \bdot b)\mdot c +(a \mdot b)\bdot c \in F^3 + A \bdot F^3 = F^3 $, hence $F^3 = I \cap F$.  Thus we have $ L \notin I$, so $I$ is proper.
\myqed

\addtocounter{mypcount}{1}\begin{proposition}\mylbl{simplephy-prop-6}
Let $R$ be a simple physical conformal superalgebra.  Then the map
$$\iota:  \Cl(V,(\cdot,\cdot)) \to \RC,$$
$$ v_1 v_2 \cdots v_r \mapsto (v_1 \mdot + v_1 \bdot)(v_2 \mdot + v_2 \bdot)\cdots (v_r \mdot + v_r \bdot) L,$$
is surjective unless $V\mdot V\mdot V =0$ with $V \ne \{0\}$.
\end{proposition}
\pr
Suppose $V=\{0\}$ and the map $\iota$ is neither zero nor surjective. Then the subspace $A \oplus F$ is closed under the $\PP{0}$ product and the $\PP{1}$ product,  so $A+F$ generates a proper ideal.
Otherwise suppose $V\mdot V\mdot V \ne \{0\}$.
Take $\IC$ the ideal of the reduced subspace $\RC$ generated by $S=V\mdot V\mdot V$.  We have $A \bdot S \subset S$ by (H3).
$ V \mdot V \bdot S \subset A\bdot S + S \subset S$ because of (H4).
Apply Proposition \ref{primalg-prop-9} taking an order so that $ x < y$ if $\wt(x) > \wt(y)$.
 Then we have $\IC\cap F = S$, so $S=F$ because $R$ is simple and $S\ne\{0\}$.
  Take a unit vector $e$ of $V$.  
$\RC$ is a $\Cl(V,(\cdot,\cdot))$-module and the Clifford action of $e$ yields an isomorphism between $\C L \oplus A$ and $V\oplus F$.   Since $V\oplus F = V\oplus S\subset \Img\iota$, so $\C L \oplus A\subset \Img\iota$, thus the map $\iota$ is surjective.
\myqed

We shall denote the conformal sub-superalgebra generated by $\Img\iota$ by $R_\iota$.

\section{ Invariants }

Let $R$ be a physical conformal superalgebra, $\RC$ the reduced subspace, $V=\RC^\hf{3}$, $A=\RC^1$ and $F=\RC^\hf{1}$.
Consider the trilinear map defined by $$ \eta: V\times V \times V \to V$$
$$(u,v,w) \mapsto u \bdot v \mdot w ,$$
and the bilinear form $(\cdot,\cdot)_{V\wedge V}$ on $V\wedge V$ defined by $(u\wedge v , w \wedge z)_{V \wedge V}L = u \bdot \eta(v, w ,z) $. 
The form is well-defined on $V \wedge V$ because
$ u \bdot v \bdot w \mdot z =  - u \bdot v \bdot z \mdot w $
and
\begin{eqnarray}
 u \bdot v \bdot w \mdot z &=& (u \bdot v)\bdot w \mdot z - v \bdot u \bdot w \mdot z\nonumber\\
 &=& - v \bdot u \bdot w \mdot z.
\end{eqnarray}
The form $(\cdot,\cdot)_{V\wedge V}$ is symmetric because
\begin{eqnarray}
 u \bdot v \bdot w \mdot z &=& - u \bdot w \bdot v \mdot z - u \bdot v \mdot w \bdot z  - u \bdot w \mdot v \bdot z + 2 u \bdot ( v \bdot w ) \mdot z\nonumber \\
&=& - w \bdot u \bdot z \mdot v - (u,v)(w,z)L - (u,w)(v,z)L + 2 (u,z)(v,w)L\nonumber\\
& & - 2 (w,v)(u,z)L  - (u,v)(w,z)L - (u,w)(v,z)L + 2 (u,z)(v,w)L\nonumber\\
&=& w \bdot z \bdot u \mdot v. 
\end{eqnarray}

The form $(\cdot,\cdot)_{V\wedge V}$ is invariant under isomorphisms of physical conformal superalgebras.

\addtocounter{mypcount}{1}\begin{proposition}\mylbl{inv-inv}
Let $f:R\to R'$ be an isomorphism of physical conformal superalgebras. $f$ induces an isometric transformation $f\wedge f:(V\wedge V, (\cdot,\cdot)_{V\wedge V}) \to (V'\wedge V', (\cdot,\cdot)_{V'\wedge V'})$ where $V' = \RC'^{\hf{3}}$.
\end{proposition}
\pr
For an isomorphism $f$ we have
\begin{eqnarray}
(u \wedge v , w \wedge z ) f(L) &=& f((u \wedge v , w \wedge z )L)\nonumber\\
 &=& f( u \bdot v \bdot w \mdot z )\nonumber\\
 &=&  f(u) \bdot f(v) \bdot f(w) \mdot f(z)\nonumber\\
 &=& ( f(u) \wedge f(v) , f(w) \wedge f(z) )f(L).
\end{eqnarray}
\myqed

Furthermore we can reconstruct the products $\mdot$ and $\bdot$ on $\iota(\Cl(V,(\cdot,\cdot)))$ from the form $(\cdot,\cdot)_{V\wedge V}$ for simple physical conformal superalgebras.

\addtocounter{mypcount}{1}\begin{proposition}\mylbl{invariants-prop-5}
Let $R$ be a physical conformal superalgebra and $\RC$ the reduced subspace.  Consider the map $\iota$ of Proposition \ref{simplephy-prop-6} and the map $\eta$.
Then  $x_{(n)}y$ is uniquely determined by the pair $( \iota , \eta)$ for all $n\in \N$ where $x, y \in R_\iota$.
\end{proposition}
\pr
The actions $V\mdot$ and $V\bdot$ are uniquely determined by $\iota$ and $\eta$, since
\begin{eqnarray}
 u \mdot v &=& \hf{1} \iota( uv - vu),\mylbl{decomp-base1}\\
 u \mdot v \mdot w &=& \iota(uvw) - \iota(\eta(u,v,w)) - (v,w)_V \iota(u),\mylbl{decomp-base2}\\
 u \bdot v \mdot w \mdot z &=& \iota(uvwz) - \iota(u\eta(v,w,z)) - (w,z)_V \iota(uv),\mylbl{decomp-base3}\\
 u \bdot v &=& \hf{1} \iota( uv + vu),\\
 u \bdot v \mdot w &=& \iota(\eta(u,v,w)),\\
 u \mdot v \bdot w \mdot x \mdot y &=& ( u \mdot v ) \bdot w \mdot x \mdot y + ( u \bdot v) \mdot w \mdot x \mdot y\nonumber\\
 &=& - \eta(w,u,v)\mdot x \mdot y  -  w \mdot \eta(x,u,v) \mdot y - w \mdot x \mdot \eta(y,u,v)\nonumber\\
 &=& -\iota(\eta(w,u,v)xy)  - \iota(w\eta(x,u,v)y) - \iota(wx\eta(y,u,v))\nonumber \\
 & &+ \iota(\eta(\eta(w,u,v),x,y))+ \iota(\eta(w,\eta(x,u,v),y)) + \iota(\eta(w,x,\eta(y,u,v)))\nonumber\\
 & &\!\!\!\!\!\!\!\!\!\!\!\!+ (x,y)_V \iota(\eta(w,u,v)) + (\eta(x,u,v),y)_V \iota(w) + (x,\eta(y,u,v))_V \iota(w),\\
 u \bdot v \bdot w \mdot x \mdot y &=& \iota(uvwxy) - \iota(uv\eta(w,x,y)) - (x,y)_V \iota(uvw)\nonumber\\
& &\!\!\!\!\!\!\!\!\!\!\!\!+\iota(\eta(w,u,v)xy)  + \iota(w\eta(x,u,v)y) + \iota(wx\eta(y,u,v)) \nonumber\\
& &\!\!\!\!\!\!\!\!\!\!\!\!- \iota(\eta(\eta(w,u,v),x,y))- \iota(\eta(w,\eta(x,u,v),y)) - \iota(\eta(w,x,\eta(y,u,v)))\nonumber\\
& &\!\!\!\!\!\!\!\!\!\!\!\!- (x,y)_V \iota(\eta(w,u,v)) - (\eta(x,u,v),y)_V \iota(w) - (x,\eta(y,u,v))_V \iota(w).
\end{eqnarray}
By Lemma \ref{primalg-lem-solve} we have $\pp{(\pp{b}{q}{c})}{p}{a}=\sum_j {r_{j}\pp{b}{q-j}{\ppp{c}{p+j}{a}}+s_{j}\pp{c}{p+q-j}{\ppp{b}{j}{a}}}$ for some $r_{j},s_{j}\in K$, hence all $(v \odot x)\odot$ are written in some $v\odot$s and  $x\odot$s where $v \in V$, $x \in \Img\iota$ and $\odot$ denotes any of $\mdot$ and $\bdot$. Since $\Img\iota = \K L \oplus V \oplus ((V \mdot V) + (V \bdot V \mdot V \mdot V)) \oplus (V \mdot V \mdot V)$, we have the results.
\myqed

Denote $\Cl^n(V,(\cdot,\cdot)) = \Span\{ v_1v_2\cdots v_k| v_i \in V,\ k\le n\}$.
(\ref{decomp-base1}), (\ref{decomp-base2}) and (\ref{decomp-base3}) imply the following.

\addtocounter{mypcount}{1}\begin{proposition}\mylbl{inv-prop-a}
Let $R$ be a simple physical conformal superalgebra and $\RC$ the reduced subspace.  For the map $\iota$ of Proposition \ref{simplephy-prop-6}, we have $\check{R_\iota} = \iota(\Cl^4(V,(\cdot,\cdot)))$.  Furthermore if $V\mdot V \mdot V = \{0\}$ then $\check{R_\iota} = \iota(\Cl^2(V,(\cdot,\cdot)))$.
\end{proposition}

\section{Classification of simple physical conformal superalgebras}

We start classification of simple physical conformal superalgebras.    We will follow the notations $\Vir$, $K_1$, $K_2$, $K_3$, $S_2$, $W_2$, and $\CKS$ given in \cite{ck} and \cite{kac0}. By the results of the preceding sections all that we have to do is listing up the left $\Cl(V,(\cdot,\cdot))$-submodules of $\Cl(V,(\cdot,\cdot))$ and the symmetric forms on $V\wedge V$ appropriate to reconstruct simple physical conformal superalgebras.

Fix a vector space $V$ with the nondegenerate inner product $(\cdot,\cdot)$ and consider an orthonormal basis $\{e_1,e_2,\cdots,e_N\}$ of $V$.  Set $D^0_k = D_k = \frac{1}{\sqrt{2}}(e_{2k-1} +  i e_{2k})$,
$D^1_k = \bar{D}_k = D_{\bar{k}} = \frac{1}{\sqrt{2}}(e_{2k-1} -  i e_{2k})$, and $D^w = D^{w_1}_1D^{w_2}_2\cdots D^{w_n}_n$ where $n=\lfloor N/2 \rfloor$, $w\in (\ZZ)^n$ and $w_i$ denotes the $i$th binary digit of $w$.
We have the following theorem for the decomposition of left $\Cl(V)$-module $\Cl(V)$.(\cite{chevalley})

\addtocounter{mypcount}{1}\begin{theorem}\mylbl{invariants-th-6}
The left $\Cl(V)$-module $\Cl(V)$ is completely reducible.
The irreducible decomposition is given as follows.
If $N=2n$ then
\begin{equation}
\Cl(V) = \bigoplus_{w \in (\ZZ)^n}{ M(w)},
\end{equation}
where $M(w) = \Cl(V)D^w$.
If $N=2n+1$ then
\begin{equation}
\Cl(V) = \bigoplus_{w \in (\ZZ)^n} {(M^+(w) \oplus M^-(w))},
\end{equation}
where $M^\pm(w) = \Cl(V)D^w (1\pm e_{N})$.
\end{theorem}

\addtocounter{mypcount}{1}\begin{proposition}
Let $R$ be a simple physical conformal superalgebra with $\dim V \le 3$.  Then $R$ is isomorphic to one of $\Vir$, $K_1$, $K_2$, $K_3$.
\end{proposition}
\pr
If $\dim V =0$ then the map $\iota$ is surjective, so $R$ is isomorphic to $\Vir$.  
Otherwise by Proposition \ref{simplephy-prop-4} $\dim F \le \left ( \matrix{{\dim V} \cr 3} \right )$ for a simple physical conformal superalgebra $R$, so $R=R_\iota$ for $V \ne \{0\}$ unless $\dim V = 3$ with $V\mdot V\mdot V=\{0\}$.
If $\dim V =1$ then we have $V\wedge V = \{0\}$, so the conformal superalgebra $R=R_\iota$ is unique, which is $K_1$.  If $\dim V =2$ then $\dim V\wedge V =1$.
Since $\bar{D}_1 \bdot D_1 \bdot \bar{D}_1 \mdot D_1 = L$ by (D5), so the only possible form $(\cdot,\cdot)_{V\wedge V}$ is $ (\bar{D}_1 \wedge D_1,\bar{D}_1 \wedge D_1)_{V\wedge V} = 1 $.  Hence $\ker \iota = \{0\}$ because otherwise the form $(\cdot,\cdot)_{V\wedge V} = 0$, thus the conformal superalgebra $R=R_\iota$ with $\dim V =2$ is unique, which is $K_2$.  
If $\dim V = 3$ then 
 $\bar{D}_1 \bdot D_1 \bdot \bar{D}_1 \mdot D_1 = L$,
 ${D_1} \bdot e_3 \bdot {D_1} \mdot e_3 = 0$,
 $\bar{D}_1 \bdot e_3 \bdot \bar{D}_1 \mdot e_3 = 0$ and
 $\bar{D}_1 \bdot e_3 \bdot {D_1} \mdot e_3 = -L$ by (D5).
So the possible form $(\cdot,\cdot)_{V\wedge V}$ is uniquely determined, which is nondegenerate.  Hence $\ker \iota = \{0\}$ because otherwise the form $(\cdot,\cdot)_{V\wedge V} $ is degenerate.  Thus the conformal superalgebra $R_\iota$ with $\dim V =3$ is unique, which is $K_3$.  $V\mdot V\mdot V \ne \{0\}$ for $K_3$, so $R=K_3$.
\myqed

Consider the polynomial ring $X^n = \C[x_1,x_2,\cdots,x_n]$ of Grassmann indeterminates.  $X^n$ is decomposed into $X^n = \bigoplus_{s\in \Z^n}{X^n_s}$ by the multidegree of polynomials.
Define the action of $\Cl(V)$ on $X^n$ by  $ D_i f = \sqrt{2}x_i f$ and $\bar{D}_i f = \sqrt{2}\partial_i f$ for $f \in X^n$ where $\partial_i$ denotes $\frac{\partial}{\partial x_i}$.  
If $N=2n$ then the action $\rho:\Cl(V) \to \End(X^n)$ is an isomorphism, so for a left ideal of $\Cl(V)$ we have an isomorphism of vector spaces 
$\rho_I: \Cl(V)/I \stackrel{\sim}\to \Hom(\cap_{w\in W} \ker \rho(D^w) , X^n)$,
where $W \subset (\ZZ)^n$ and $I =\bigoplus_{w\in W} M(w)$.
Hence $\Cl(V)/I$ is decomposed into 
\begin{equation}
 \Cl(V)/I = \bigoplus_{ t \in \{-1,0,1\}^n }{(\Cl(V)/I)_t},
\end{equation}
where $(\Cl(V)/I)_t = \{ u \in (\Cl(V)/I)|\ \rho_I(u)(\cap_{w\in W}\ker \rho(D^w) \cap X^n_s) \subset X^n_{s+t}\ \mbox{for all}\ s\in \Z^n\}$.  Denote the projections $\pi^I_t: \Cl(V)/I \to (\Cl(V)/I)_t$.
If $N=2n+1$ then we have a decomposition as left $\Cl(V)$-modules 
$\Cl(V) = \Cl(V/{\C e_{N}}) (1+e_{N}) \oplus \Cl(V/{\C e_{N}}) (1-e_{N})$.
So we have an isomorphism 
$\rho_I: \Cl(V) / I \stackrel{\sim}\to \Hom(\cap_{w\in W^+} \ker \rho(D^w) , X^n) \oplus \Hom(\cap_{w\in W^-} \ker \rho(D^w) , X^n)$,
 where $ W^\pm \subset (\ZZ)^n $ and $I =( \bigoplus_{w\in W^+} M^+(w) ) \oplus (\bigoplus_{w\in W^-} M^-(w) )$.
Hence the decomposition is
\begin{equation}
 \Cl(V)/I = \bigoplus_{ t \in \{-1,0,1\}^n }{(\Cl(V)/I)^+_t} \oplus \bigoplus_{ t \in \{-1,0,1\}^n }{(\Cl(V)/I)^-_t}.
\end{equation}
Denote the projections $\pi^{I\,\pm}_t: \Cl(V)/I \to (\Cl(V)/I)^\pm_t$.

Consider $\alpha_{i,j}\in \C$ defined by $\alpha_{i,j}L = D_i\bdot\bar{D}_i\bdot{D_j}\mdot\bar{D}_j$ for $i,j \in \{1,2,\cdots,n\}$ and $\beta_{i,j,k,l}\in \C$ by $\beta_{i,j,k,l}L = D_i\bdot{D_j}\bdot{D_k}\mdot{D_l}$ for $i,j,k,l \in \{1,\bar{1},2,\bar{2},\cdots,n,\bar{n}\}$.  We have\\ $D_j \bdot D_i \bdot \bar{D}_i \mdot \bar{D}_j = (1+\alpha_{i,j})L$, $D_j \bdot \bar{D}_i \bdot D_i \mdot \bar{D}_j = (1-\alpha_{i,j})L$ and $\alpha_{i,j} = \alpha_{j,i}$.

\addtocounter{mypcount}{1}\begin{proposition}
Let $R$ be a simple physical conformal superalgebra with $\dim V \ge 4$.  Then $\dim V$ is one of $4$, $6$, $8$.  Furthermore if $V\mdot V\mdot V =\{0\}$ then $\dim V =4$.
\end{proposition}
\pr
Suppose given a simple physical conformal superalgebra $R$ with $\dim V = 2n+1$ where $n \ge 2$.
  For an arbitrary $u\in \ZZ$ we have
\begin{eqnarray}
0&=&-(e_{2n+1}\mdot\bar{D}_i)\bdot(D_i\mdot D^u_j) -(D_i\mdot D^u_j)\bdot (e_{2n+1}\mdot\bar{D}_i)\nonumber\\
&=&(D_i\bdot e_{2n+1} \mdot \bar{D}_i)\mdot D^u_j + D_i\mdot (D^u_j \bdot e_{2n+1} \mdot \bar{D}_i)\nonumber\\
& & + (e_{2n+1}\bdot D_i \mdot {D^u_j})\mdot \bar{D}_i + e_{2n+1}\mdot (\bar{D}_i \bdot D_i \mdot {D^u_j}).
\end{eqnarray}
Apply $\pi^{\ker\iota\,+}_t + \pi^{\ker\iota\,-}_t$ to the both sides where all digits of $t$ are $0$ except for $t_j=(-1)^u$.  Then we have $((-1)^u\alpha_{i,j}+2) e_{2n+1}\mdot D^u_j = 0$.  If $e_{2n+1}\mdot D^u_j = 0$ then $ 0 = e_{2n+1} \bdot e_{2n+1}\mdot D^u_j = D^u_j$, so we have $\alpha_{i,j}= -(-1)^u 2$ for an arbitrary $u\in\ZZ$ where $i,j=1,2,\cdots,n$. 
Hence $R$ does not exist.

Suppose given a simple physical conformal superalgebra with $\dim V =2n$ and $n > 4$.
Consider $S\subset (\ZZ)^n$ such that $\ker \iota=I=\bigoplus_{s\in S}M(s)$.
By Proposition \ref{inv-prop-a} $\Cl(V)/I = \Cl^4(V)/I$, so we have $\pi^I_t = 0$ if $\#\{k\in\N|\ t_k\ne0\} > 4$.
Hence $\iota(D^w) = \pi^I_t(\iota(D^w)) = 0$ for an arbitrary $w\in(\ZZ)^n$ where $t_i=(-1)^{w_i}$, so $S = (\ZZ)^n$, that is, the map $\iota$ is the $0$ map, thus we have the result.

If $V\mdot V\mdot V =\{0\}$ then we have $\dim V = 4$ in the same way by $\Cl(V)/I = \Cl^2(V)/I$.
\myqed

\addtocounter{mypcount}{1}\begin{proposition}\mylbl{inv-prop-y}
A simple physical conformal superalgebra $R$ with $\dim V =6$ is isomorphic to $\CKS$
\end{proposition}
\pr
Suppose given a simple physical conformal superalgebra $R$ with $\dim V = 6$.
The map $\iota$ is surjective.
For all $i,j,k\in \N$ satisfying $\{i,j,k\}=  \{1,2,3\}$ we have
\begin{eqnarray}
0 &=& (\Eb{i}{j})  \bdot (\F{k}{i}) + (\F{k}{i}) \bdot (\Eb{i}{j})\nonumber \\
&=& (D_i \bdot \bar{D}_i\mdot \bar{D}_j)\mdot \bar{D}_k+ D_i\mdot(\bar{D}_k\bdot\bar{D}_i\mdot\bar{D}_j)\nonumber \\
& & +(\bar{D}_i \bdot {D}_i\mdot \bar{D}_k)\mdot \bar{D}_j+ \bar{D}_i\mdot(\bar{D}_j\bdot{D}_i\mdot\bar{D}_k)\nonumber\\
&=& (\alpha_{i,j} + \alpha_{i,k})\bar{D}_j\mdot\bar{D}_k + \cdots.
\end{eqnarray}
Applying $\pi^{\ker \iota}_t$ to the both sides where all trial digit of $t$ is $0$ except for the $i$th and the $j$th, $t_i=t_j=-1$, we have $(\alpha_{i,j} + \alpha_{i,k})\bar{D}_j\mdot\bar{D}_k = 0$. Hence we have $(\alpha_{i,j} + \alpha_{i,k})(\alpha_{j,k} +1) = 0$ applying ${D_j}\bdot{D_k}\bdot$ to the both sides.
Similarly for $ 0 = (\F{i}{j}) \bdot (\F{k}{i}) + (\F{k}{i}) \bdot (\F{i}{j}) $ and the term $D_j\mdot \bar{D}_k$ we have $(\alpha_{i,j} - \alpha_{i,k})(\alpha_{j,k} -1) = 0$.
It is easily checked that the solutions of the above equations are $\ba =(\alpha_{1,2},\alpha_{2,3},\alpha_{3,1}) = (0,0,0)$, $\ba=(-1,-1,-1)$, $\ba=(-1,1,1)$, $\ba=(1,-1,1)$, $\ba=(1,1,-1)$.
If $\ba = (-1,-1,-1)$ then we have $D_i\mdot D_j = 0$ and $\bar{D}_i\mdot\bar{D}_j = 0$ for all $i,j = 1,2,3$, which implies $D^{w_1}_i\mdot D^{w_2}_j\mdot D^{w_3}_k = 0$ for all $w = (w_1,w_2,w_3) \in (\ZZ)^3$ where $\{i,j,k\}=\{1,2,3\}$, so $R =0$.  In the same way $R=0$ for each $\ba=(-1,1,1)$, $\ba=(1,-1,1)$, $\ba=(1,1,-1)$.
If  $\ba = (0,0,0)$ then the set $\{D^{w_1}_i\mdot D^{w_2}_j | i\ne j,\ w\in(\ZZ)^2\}$ is linearly independent.  Considering the relations $ 0 = (D^{w_1}_i\mdot D^{w_2}_j) \bdot (D^{w'_1}_k \mdot D^{w'_2}_l) +  (D^{w'_k}_i\mdot D^{w'_2}_k) \bdot (D^{w_1}_i \mdot D^{w_2}_j) $ for all $w,w'\in (\ZZ)^2$ and $i,j,k,l \in\{1,2,3\}$ in the same way, we have $\beta_{i,j,k,l} = 0$ if $\{i,j,k,l\}\ne \{s,\bar{s},t,\bar{t}\}$ for any $s,t\in \{1,2,3\}$.   Consider the map $\iota$ for this case.  By Proposition \ref{inv-prop-a} we may assume $\bar{D}_1\bar{D}_2\bar{D}_3 \in \ker \iota$.  If we have ${D^{w_1}_1}\mdot{D^{w_2}_2}\mdot{D^{w_3}_3} \ne 0$ for some $w \in\{ (1,0,0),(0,1,0),(0,0,1)\}$ then $R$ is not simple by Proposition \ref{simplephy-prop-4}.  Hence ${D^{w_1}_1}\mdot{D^{w_2}_2}\mdot{D^{w_3}_3} = 0$ for all $w \in\{ (1,1,1),(1,0,0),(0,1,0),(0,0,1)\}$.  Thus the simple physical conformal superalgebra structure on this space is uniquely determined, which is $\CKS$.
\myqed

\addtocounter{mypcount}{1}\begin{proposition}
Simple physical conformal superalgebras with $\dim V = 8$ do not exist.
\end{proposition}
\pr
Suppose given a simple physical conformal superalgebra $R$ with $\dim V = 8$.
The map $\iota$ is surjective.
Then we have $(\alpha_{i,j}+\alpha_{j,k})(\alpha_{i,k}+1)=0$ and $(\alpha_{i,j}-\alpha_{j,k})(\alpha_{i,k}-1)=0$ for all distinct $i,j,k$.  It is easily checked that the set of solutions $\ba = (\alpha_{1,2},\alpha_{1,3},\alpha_{1,4},\alpha_{2,3},\alpha_{2,4},\alpha_{3,4})$ is $\{(0,0,0,0,0,0), (1,1,1,-1,-1,-1), (1,1,-1,-1,1,1), (1,-1,1,1,-1,1), (1,-1,-1,1,1,1),\\ (-1,1,-1,1,-1,1), (-1,1,1,1,1,-1), (-1,-1,1,-1,1,1), (-1,-1,-1,-1,-1,-1)\}$.  For the non-zero solutions we have $D^{w_1}_i D^{w_2}_j D^{w_3}_k D^{w_4}_l = 0$ for all $w\in (\ZZ)^4$ where $\{i,j,k,l\}=\{1,2,3,4\}$, hence $R = 0$.
If $\ba =0$ then the set $\{D^{w_1}_i\mdot D^{w_2}_j | i\ne j,\ w\in(\ZZ)^2\}$ is linearly independent, so $\beta_{i,j,k,l} = 0$ if $\{i,j,k,l\}\ne \{s,\bar{s},t,\bar{t}\}$ for any $s,t\in\{1,2,3,4\}$, which implies  $\iota(D^{w_1}_i D^{w_2}_j D^{w_3}_k D^{w_4}_l) = 0$ if and only if $D^{w_1}_i\bdot D^{w_2}_j\mdot D^{w_3}_k\mdot D^{w_4}_l = 0$.
Suppose $\iota(D^{w_1}_i D^{w_2}_j D^{w_3}_k D^{w_4}_l) = 0$ for some $w \in (\ZZ)^4$ where $\{i,j,k,l\}=\{1,2,3,4\}$.  Then 
$ 0 = D^{w_1+1}_i \mdot D^{w_1}_i\bdot D^{w_2}_j\mdot D^{w_3}_k\mdot D^{w_4}_l
= D^{w_2}_j\mdot D^{w_3}_k\mdot D^{w_4}_l$
because of (H4), so we have $\iota(D^{w_1+1}_i D^{w_2}_j D^{w_3}_k D^{w_4}_l) = 0$.
Since $\ker \iota \ne \{0\}$ by Proposition \ref{inv-prop-a}, thus $R=0$.
\myqed

\addtocounter{mypcount}{1}\begin{proposition}\mylbl{invariants-prop-7}
For a simple physical conformal superalgebra with $\dim V=4$ the form $(\cdot,\cdot)_{V\wedge V}$ is as given on table \ref{n4table} for some $\alpha \in \C$.
\begin{table}[htbp]
\caption{The form $(\cdot,\cdot)_{V\wedge V}$ for $\dim V =4$.}
\mylbl{n4table}
\begin{center}
\begin{tabular}{c|cccccc}
$(\cdot,\cdot)_{V\wedge V}^\alpha$ & $\bar{D}_1\wedge D_1$ & $\bar{D}_2\wedge D_2$ & $D_1\wedge D_2$ & $\bar{D}_1\wedge D_2$ & $D_1\wedge \bar{D}_2$ & $\bar{D}_1\wedge \bar{D}_2$ \\ \hline
$\bar{D}_1\wedge D_1$ & $1$ & $\alpha$ & $0$ & $0$ & $0$ & $0$ \\
$\bar{D}_2\wedge D_2$ & $\alpha$ & $1$ & $0$ & $0$ & $0$ & $0$ \\
$D_1\wedge D_2$ & $0$ & $0$ & $0$ & $0$ & $0$ & $-(1+\alpha)$ \\
$\bar{D}_1\wedge D_2$ & $0$ & $0$ & $0$ & $0$ & $-(1-\alpha)$ & $0$ \\
$D_1\wedge \bar{D}_2$ & $0$ & $0$ & $0$ & $-(1-\alpha)$ & $0$ & $0$ \\
$\bar{D}_1\wedge \bar{D}_2$ & $0$ & $0$ & $-(1+\alpha)$ & $0$ & $0$ & $0$ \\
\end{tabular}
\end{center}
\end{table}
\end{proposition}
\pr
The table \ref{n4table} is obtained by using the following formulae:
\begin{eqnarray}
D^{a}_1 \bdot D^{a}_1 \mdot D^{b}_2 &=& 0,\\
D^{a}_2 \bdot D^{a}_2 \mdot D^{b}_1 &=& 0,\\
D^{a}_1 \bdot D^{a+1}_1 \mdot D^{b}_2 &=& - D^{a}_1 \bdot D^{b}_2 \mdot D^{a+1}_1\nonumber\\
 &=& D^{a}_1 \mdot D^{b}_2 \bdot D^{a+1}_1 + D^{b}_2 \bdot D^{a}_1 \mdot D^{a+1}_1\nonumber \\
& & \quad + D^{b}_2 \mdot D^{a}_1 \bdot D^{a+1}_1 - 2(D^{a}_1 \bdot D^{b}_2) \mdot D^{a+1}_1\nonumber\\
&=& D^{b}_2 + (-1)^{a+b} \alpha D^{b}_2, \\
D^{a}_2 \bdot D^{a+1}_2 \mdot D^{b}_1 &=& D^{b}_1 + (-1)^{a+b} \alpha D^{b}_1, 
\end{eqnarray}
where $ a,b \in \ZZ $
\myqed 
\addtocounter{mypcount}{1}\begin{remark}{\rm \mylbl{invariants-rem-8}
The form $(\cdot,\cdot)_{V \wedge V}$ of Proposition \ref{invariants-prop-7} is given by
\begin{equation}
(e_i\wedge e_j , e_k\wedge e_l)_{V\wedge V} = -\alpha \epsilon_{ijkl} + \delta_{jk}\delta_{il}- \delta_{ik}\delta_{jl}, \label{alpha_onb}
\end{equation}
where $\epsilon_{ijkl}$ is antisymmetric with $\epsilon_{1234}=1$.
}\end{remark}

On the other hand we have the following proposition.

\addtocounter{mypcount}{1}\begin{proposition}\mylbl{invariants-prop-nfa}
A physical conformal superalgebra structure exists on $\Cl(V)$ where $\dim V = 4$ with the form $(\cdot,\cdot)_{V\wedge V}$ described in Table \ref{n4table} for an each $\alpha\in\C$.
\end{proposition}
\pr
By Proposition \ref{invariants-prop-5} the form $(\cdot,\cdot)_{V\wedge V}$ determines the products $\mdot$ and $\bdot $ on $\Cl(V)$ for an arbitrary $\alpha\in\C$.
It is easily checked that they have all properties (H0-6).  By Proposition \ref{dotalg-prop-8.1a} a conformal superalgebra structure is determined on $\Cl(V)$ for each $\alpha \in\C$.
\myqed

We shall denote thus obtained family of physical conformal superalgebras by $\{\NFA\}_{\alpha\in\C}$.  
$\NFA$ is equivalent to $N_4^0$ changing the conformal vector $L$ to ${L_\alpha} =  L - \frac{\alpha}{2}\partial e_1\bdot e_2\mdot e_3\mdot e_4$ except for $\alpha^2 =1$.  For $\alpha^2 = 1$ we shall denote $N_4 = (N_4^0, (n), {L_1})$, which is isomorphic to $(N_4^0, (n), L_{-1})$.

\addtocounter{mypcount}{1}\begin{note}\mylbl{invariants-note-k4sub}{\rm
The conformal superalgebra $K_4$ is written down in \cite{kac2}.  The physical conformal superalgebra $N_4^0$ is isomorphic to the subalgebra of $K_4$ generated by the primary vectors other than $\xi_1\xi_2\xi_3\xi_4$.
}\end{note}

\addtocounter{mypcount}{1}\begin{proposition}\mylbl{invariants-prop-12}
$\NFA$ and $\NFB$ are isomorphic if and only if $\alpha^2= \beta^2$.
\end{proposition}
\pr
Set $E_1=\bar{D}_1\wedge D_1$, $E_2=\bar{D}_2\wedge D_2$, $E_3=D_1\wedge D_2$, $E_4=\bar{D}_1\wedge D_2$, $E_5=D_1\wedge \bar{D}_2$, $E_6=\bar{D}_1\wedge \bar{D}_2$. 
The characteristic polynomial of the matrix $M_{i,j}=(E_i,E_j)_{V\wedge V}$ is $ ((t+1)-\alpha^2)((t-1)^2-\alpha^2)^2$, which is invariant under automorphisms by Proposition \ref{inv-inv}.  Hence if $\NFA$ and $\NFB$ are isomorphic then $\alpha^2 = \beta^2$.  Conversely suppose $\beta = -\alpha$.  Consider the map $f: V \to V$ defined by $f(e_1) =e_2 $, $f(e_2) =e_1$, $f(e_k) =f(e_k)$ for all $k>2$.  Because $\ker \iota = \{0\}$, $f$ extends to an automorphism of conformal superalgebra, which maps $\alpha$ to $-\alpha$.
\myqed

For their simplicity we have the following.

\addtocounter{mypcount}{1}\begin{proposition}\mylbl{inv-prop-s}
$\NFA$ is simple if and only if $\alpha^2 \ne 1$.
\end{proposition}
\pr
 Consider the bilinear form $\langle\cdot,\cdot\rangle$ on $V\wedge V\wedge V$ defined by $\langle u_1\wedge u_2 \wedge u_3, v_1\wedge v_2 \wedge v_3 \rangle L= u_3\bdot u_2 \bdot u_1 \bdot v_1 \mdot v_2 \mdot v_3$.   Denote $f_{i,j}=D_i \wedge \bar{D}_i\wedge D_j$ and $\bar{f}_{i,j}=\ \bar{D}_i\wedge D_i \wedge \bar{D}_j$ and take the basis
$\{f_{i,j},\bar{f}_{i,j}|\ \{i,j\} = \{1,2\} \}$ of $V\wedge V \wedge V$.  Then we have $\langle \bar{f}_{i,j}, f_{k,l} \rangle = \langle f_{i,j}, \bar{f}_{k,l} \rangle = (1-\alpha^2) \delta_{i,k}\delta_{j,l}$, $\langle f_{i,j}, f_{k,l} \rangle = \langle \bar{f}_{i,j}, \bar{f}_{k,l} \rangle = 0 $.  So the form $\langle \cdot,\cdot \rangle$ is symmetric, and is nondegenerate if and only if $\alpha^2 \ne1$.  By Proposition \ref{simplephy-prop-4}, $\NFA$ is simple if and only if $\alpha^2 \ne 1$.
\myqed

In particular $N_4^0$ is simple, so we have the following corollary.

\addtocounter{mypcount}{1}\begin{corollary}\mylbl{inv-cor-n4}
$N_4$ is simple.
\end{corollary}

\addtocounter{mypcount}{1}\begin{note}{\rm 
A one-parameter family of superconformal algebras that is called the large $N=4$ superconformal algebra is written down in \cite{STP}.  
In (2), (3), (4) of \cite{STP} set $\gamma = (\beta +1)/2$ and replace the central terms by $0$.  Fix the conformal vector $L(z) = \sum_{n\in\N}{L_{n}z^{-n-2}}$.
If $\beta^2 \ne 1$ then the centerless large $N=4$ superconformal algebra is isomorphic to $\NFB$ by
\begin{eqnarray}
G_a &=& \sqrt{2}e_a,\nonumber\\
A^{\pm\,1} &=& \frac{-1}{2(1\pm\beta)} \left( e_2\mdot e_3 \pm e_1\mdot e_4\right),\nonumber\\
A^{\pm\,2} &=& \frac{1}{2(1\pm\beta)} \left( e_1\mdot e_3 \mp e_2\mdot e_4\right),\nonumber\\
A^{\pm\,3} &=& \frac{-1}{2(1\pm\beta)} \left( e_1\mdot e_2 \pm e_3\mdot e_4\right),\\
Q^a &=& \frac{1}{\sqrt{2}(1-\beta^2)}(-1)^a e_{b_1}\mdot e_{b_2} \mdot e_{b_3},\nonumber\\
U &=& \frac{-1}{1-\beta^2} e_1\bdot e_2 \mdot e_3 \mdot e_4,\nonumber
\end{eqnarray}
where $\{a,b_1,b_2,b_3\} = \{1,2,3,4\}$ and $b_1 < b_2 < b_3$.
For $\beta=\pm 1$ the large $N=4$ superconformal algebra is isomorphic to $N_4$.
}\end{note}

\addtocounter{mypcount}{1}\begin{note}{\rm 
 The action of the Lie algebra $(A,(0))$ on $V$ is not faithful for $\NFA$ since $e_1 \bdot e_2 \mdot e_3 \mdot e_4$ acts on $V$ trivially. 
The ideal generated by $e_1\bdot e_2 \mdot e_3 \mdot e_4$ is  $R$ itself.
As is discussed in (4.12) of \cite{kac2}, for a unit vector $u\in V$ one has an $A_u$-module isomorphism $ u\mdot{}: A_u \stackrel{\sim}\to F$ with the inverse map $u\bdot{}: F \stackrel{\sim}\to A_u$ where $A_u = \{ a\in A|\ a\bdot u=0 \}$ in our terminology.  $A_{e_1}$ is spanned by  $A^{+\,1}+A^{-\,1}$, $A^{+\,2}-A^{-\,2}$, $A^{+\,3}-A^{-\,3}$ and $U$, so the condition that $F$ is isomorphic to $A_u$ as $A_u$-modules is also satisfied here.
}\end{note}

\addtocounter{mypcount}{1}\begin{proposition}\mylbl{invariants-th-14}
A simple physical conformal superalgebra $R$ with $\dim V = 4$ is isomorphic to one of $S_2$, $W_2$, $N_4$ and $\NFA$ for some $\alpha \in \C$ where $\alpha \in (\C/\{\pm 1\}) \setminus \{[1]\}$.
\end{proposition}
\pr
Since the Clifford action of a unit vector in $V$ yields an isomorphism between the even subspace of $\Img\iota$ and the odd subspace of $\Img\iota$, we have $\dim \Img\iota \ge 2 \dim V$, so $\dim \Img\iota$ is one of $8$, $12$ and $16$.

If $\dim \Img\iota=16$ then the map $\iota$ is injective, so $R_\iota$ is isomorphic to $\NFA$ for some $\alpha\in \C$.  So $\iota$ is surjective because $V\mdot V\mdot V \ne \{0\}$ for all $\NFA$s.  $\NFA$ is simple if and only if $\alpha^2\ne 1$, hence $R$ is isomorphic to $\NFA$ for some $\alpha^2 \ne 1$.

If $\dim \Img\iota = 12$ then $\dim \ker \iota = 4$, so
we may assume $\ker \iota = M(00)$, which implies $D^0_1 \mdot D^0_2 = 0$.  $\alpha = -1$ by Proposition \ref{invariants-prop-7}, hence $R_\iota$ is uniquely determined, which is neither simple nor with $V\mdot V\mdot V=\{0\}$.  So simple physical conformal superalgebras with $\dim \Img\iota = 12$ do not exist.

 If $\dim \Img\iota = 8$ then $\dim \ker \iota = 8$.
We may assume $M(00) \subset \ker\iota$, which implies $D^0_1 \mdot D^0_2 = 0$.  $\alpha = -1$ by Proposition \ref{invariants-prop-7}, so $\ker\iota =M(00) \oplus M(11)$.  Hence $R_\iota$ is uniquely determined, which is $S_2$.
In particular the Lie algebra $(V\mdot V, \bdot)$ and its action on $V$ by the product $\bdot$ is uniquely determined.
Consider the pairing $J: V\wedge V\wedge V \times F \to \C$ defined by $J(v^1\wedge v^2\wedge v^3 , f)L = v^1\bdot v^2 \bdot v^3 \bdot f$.  
By (D6) we have $J(a\cdot \omega, f) + J(\omega , a\cdot f) = 0$ for all $a\in A$, $f\in F$ and $\omega\in V\wedge V\wedge V$  where $V$ and $F$ are supposed to be $(A,\bdot)$-modules by the product $\bdot$ and so is $V\wedge V\wedge V$ by derivation.
Once $J$ is determined, $u \bdot v \bdot f$ is uniquely determined for all $u, v \in V$ and $f \in F$, so the action of $A= V\mdot V + V\bdot F$ on $V$ by the product $\bdot$ is uniquely determined.
The pairing $J$ is $(A,\bdot)$-invariant and if $J(\omega,f)=0$ for all $\omega \in V\wedge V\wedge V$ then $f=0$, so the action of $A$ on $F$ by the product $\bdot$ is uniquely determined by the action of $A$ on $V\wedge V\wedge V$, which determines the product $V\mdot A$ because of (H4).
By Lemma \ref{primalg-lem-solve} we have $\pp{(\pp{b}{q}{c})}{p}{a}=\sum_j {r_{j}\pp{b}{q-j}{\ppp{c}{p+j}{a}}+s_{j}\pp{c}{p+q-j}{\ppp{b}{j}{a}}}$ for all $a, b, c \in \RC$ for some $r_{j},s_{j}\in \C$, so the simple physical conformal superalgebra structure on $R$ is uniquely determined by the pairing $J$.
Consider a $A$-submodule $J^0 = \{ \omega \in V\wedge V\wedge V |\ J(\omega, f)=0\ \mbox{for all}\ f\in F \}$.
$J^0 \ne V\wedge V\wedge V$ because $R$ is simple.
The $(V\mdot V, \bdot)$-module $V\wedge V\wedge V$ is decomposed into two $2$-dimensional irreducible modules, so $\dim J^0$ is either $0$ or $2$.
If $\dim J^0= 0$ then $R$ is isomorphic to $N_4$.
If $\dim J^0= 2$ then we can choose a basis $\{D_1, D_2, \bar{D}_1, \bar{D}_2 \}$ of $V$ that satisfies $(D_i, D_j) = (\bar{D}_i, \bar{D}_j) = \delta_{i j}$ and $(D_i, \bar{D}_j) = 0$ so that $J^0 = \Span\{D_1\wedge \bar{D}_1\wedge D_2, D_1\wedge D_2 \wedge \bar{D}_2\}$ and $\ker\iota=M(00)\oplus M(11)$, hence $R$ is unique if exists, which is $W_2$.
\myqed

Hence we have the complete list of simple physical conformal superalgebras.
\addtocounter{mypcount}{1}\begin{theorem}\mylbl{classification}
A simple physical conformal superalgebras is isomorphic to one of
$\Vir$, $K_1$, $K_2$, $K_3$, $S_2$, $W_2$, $N_4$, $\NFA$ and $\CKS$,
where $\alpha \in (\C/\{\pm 1\}) \setminus \{[1]\}$.
\end{theorem}

If conformal superalgebras $R$ and $R'$ are equivalent then the Lie superalgebras $(R/\partial R, (0) )$ and $(R'/\partial R', (0) )$ are isomorphic.
Any pair of Lie superalgebras $(\PL{\Vir}, (0))$, $(\PL{K_1}, (0))$, $(\PL{K_2}, (0))$, $(\PL{K_3}, (0))$, $(\PL{S_2}, (0))$, $(\PL{W_2}, (0))$ and $(\PL{N_4^0}, (0))$ is not isomorphic,
while $N_4$ and $\NFA$s are equivalent to $N_4^0$ except for $\alpha^2 =1$.
Hence we have the following corollary.

\addtocounter{mypcount}{1}\begin{corollary}
A simple physical conformal superalgebras is equivalent to one of
$\Vir$, $K_1$, $K_2$, $K_3$, $S_2$, $W_2$, $N_4^0$ and $\CKS$.
\end{corollary}

\end{document}